\documentclass[a4paper,12pt]{article}

\usepackage{amsmath ,amsfonts,  amssymb, amscd, amsthm} 
\usepackage{bbm}
\usepackage{mathrsfs}
\usepackage[font=small,skip=0pt]{caption}

\usepackage{csquotes, color}
\usepackage{colortbl}
\usepackage[numbers]{natbib}
\usepackage{verbatim}
\usepackage[toc,page]{appendix}
\usepackage{caption}
\usepackage{textcomp}
\usepackage[inline]{enumitem}
\usepackage{tcolorbox}
\usepackage{cancel}
\usepackage{relsize}
\usepackage{bm}

\usepackage{authblk}

\usepackage{algorithm}
\usepackage{algorithmic}
\usepackage{scrextend}
\addtokomafont{labelinglabel}{\sffamily}

\usepackage{tikz}
\usepackage{tikz-cd}
\usetikzlibrary{arrows.meta}
\usetikzlibrary{positioning}
\usetikzlibrary{shadows}
\usetikzlibrary{calc, intersections}
\usetikzlibrary{shapes.multipart}
\usetikzlibrary{patterns}

\DeclareMathOperator\erfc{erfc}

\newtheorem*{definition*}{Definition}

\newtheorem*{example*}{Example}

\newtheorem*{conditional*}{Conditional}

\newtheorem*{remark*}{Remark}

\tikzset{Myarrow/.style={very thin, arrows=-Latex}}

\definecolor{Gray}{gray}{0.85}
\definecolor{LightCyan}{rgb}{0.88,1,1}

\newcolumntype{a}{>{\columncolor{Gray}}c}
\newcolumntype{b}{>{\columncolor{Gray}}r}
\newcolumntype{d}{>{\columncolor{white}}c}

\newcommand{\cC}{{\cal C}}

\newcommand{\cG}{{\cal G}}

\newcommand{\cI}{{\cal I}}

\newcommand{\cP}{{\cal P}}

\newcommand{\cS}{{\cal S}}

\newcommand{\scrF}{{\mathscr F}} 
\newcommand{\scrG}{{\mathscr G}} 
 
\newcommand{\scrL}{{\mathscr L}}

\usepackage{hyperref}
\hypersetup{
    colorlinks=true,
    citecolor=black,
    linkcolor=black,
    filecolor=cyan,      
    urlcolor=blue,
}
 
\newenvironment{myitemize}
{ \begin{itemize}
    \setlength{\itemsep}{0pt}
    \setlength{\parskip}{0pt}
    \setlength{\parsep}{0pt}     }
{ \end{itemize}                  } 

\newenvironment{myenumerate}
{ \begin{enumerate}
    \setlength{\itemsep}{0pt}
    \setlength{\parskip}{0pt}
    \setlength{\parsep}{0pt}     }
{ \end{enumerate}                  }

\providecommand{\keywords}[1]
{
  \small	
  \textbf{\textit{Keywords---}} #1
}

\title{
Infinitely Divisible Distributions and \\ 
Commutative Diagrams
}
\author{Nomvelo Karabo Sibisi \\ {\small {\tt sbsnom005@myuct.ac.za}}}
\date{\today}

\begin{document}
\maketitle
\thispagestyle{empty}


\begin{abstract}
\noindent
We study infinitely divisible (ID) distributions on the  nonnegative  half-line $\mathbb{R}_+$.
The L\'{e}vy-Khintchine representation  of such distributions is well-known.
Our primary contribution is to cast the probabilistic objects and the relations  amongst them in 
a unified visual form  that we refer to as the  L\'{e}vy-Khintchine commutative diagram (LKCD). 
While it is  introduced as a representational tool, the LKCD facilitates the exploration of new ID distributions
and may thus also be looked upon, at least in part,  as a discovery tool.
The basic object  of the study is the gamma distribution.
Closely allied to this is the $\alpha$-stable  distribution on $\mathbb{R}_+$ for $0<\alpha<1$,
which we regard as arising from the gamma distribution rather than as a separate object.
It is indeed often characterised as an instance of a class of ID distributions known as generalised gamma convolutions (GGCs).
We make use of convolutions and mixtures of gamma and stable densities to generate densities of other GGC distributions,
with particular cases involving Bessel, confluent hypergeometric, Mittag-Leffler and parabolic cylinder functions.
We present all instances as LKCD representations.
\end{abstract}

\section{Introduction}
\label{sec:intro}
We are interested in 
infinitely divisible (ID) distributions on the  nonnegative  half-line $\mathbb{R}_+$,
that we shall occasionally refer to as `positive ID distributions' for short, even though the term is a bit  imprecise.
Such distributions find  application in numerous settings involving positive (nonnegative)  random phenomena  that are additive over 
partitions of their  domain of definition, the domain being one-dimensional (notably time) or multidimensional
(physical or more abstract space).
Given that the theory of positive ID distributions  is well-established (as briefly reviewed below), 
it is probably fair to say that worthwhile  contributions overwhelmingly lie in application.
Arguably though, there is always room for new or enhanced perspectives on theoretical foundations that are, for example,
 more visual than purely symbolic in nature.
 Such perspectives may, in turn, enrich basic understanding of the field and inform further application.
We consider this paper to lie primarily in this category of contribution,
with pointers to application  from insights that arise.

 ID distributions feature prominently in the book by Feller~\cite{Feller2}, while the one by Steutel and van Harn~\cite{SteutelvanHarn} is
 exclusively dedicated to the topic of infinite divisibility.
 Well-known examples of ID distributions  on real variables are the Gaussian and Cauchy distributions
(they also happen to be the two instances of stable distributions on $\mathbb{R}$ with known closed form densities).
  In the case of variables on $\mathbb{R}_+$, the gamma distribution is the classical example, 
 along with the  family of stable distributions on $\mathbb{R}_+$ 
 (`positive stable distributions').
 The latter two examples are related in the sense  that the gamma distribution might be described as a \textquote{founding member}
 of a  family  of distributions known as the generalised gamma convolutions (GGCs),  
 to which the positive stable distributions also belong.
 
 Thorin~\cite{Thorin77a, Thorin77b} introduced the GGC class (which belongs to the broader class of positive ID distributions)
 as he sought to prove the infinite divisibility of the Pareto  and the  log-normal distributions. 
GGC theory was subsequently studied in depth in the text by Bondesson~\cite{Bondesson}, 
exploring the powerful ramifications of the GGC concept.
The more recent survey of the GGC class by James~{\it et al.}~\cite{JamesRoynetteYor}  includes  theory and examples.
 
 All ID distributions are characterised by 
 the L\'{e}vy-Khintchine representation.
 The primary purpose of this paper is to introduce a visualisation of infinite divisibility  by 
 casting the L\'{e}vy-Khintchine representation as a commutative diagram, a construct borrowed from category theory.
  In our view, the commutative diagram lends a welcome perspective 
 to  what can often be a bewildering morass of equations in the study of ID distributions.
 
The   L\'{e}vy-Khintchine commutative diagram (LKCD) might be said to be 
 an organising principle that displays the  objects of an infinitely divisible probabilistic structure as vertices connected by arrows
 denoting relationships between objects.
 The passage from one  object to another is path-independent.
Our experience has been that  the assignment of arrow mappings such that path-independence holds can trigger a thought process 
about  ID/GGC structure that does not readily arise in the absence of the commutative diagram setting.
In some instances, this has prompted novel ideas on the representation of known densities.
This will be especially apparent when we 
describe  the fractional gamma distribution in terms of the
parabolic  cylinder function alongside the commonly used Mittag-Leffler function and its 
three-parameter generalisation known as the Prabhakar function.

Hence the LKCD can  facilitate the discovery of novel probabilistic representations.
It may not be a discovery tool in its own right,  but we believe it to be a worthy addition to the study of infinite divisibility.

\subsection{Structure of Paper}
\label{sec:structure}

Along with the well-known gamma distribution,  the stable distribution is central to the paper.
We start by introducing the stable distribution  in the conventional way as a standalone object.
We summarise infinite divisibility on $\mathbb{R}_+$ in Section~\ref{sec:ID}. 
We then  introduce the L\'{e}vy-Khintchine commutative diagram (LKCD) in Section~\ref{sec:LKCD}, 
with gamma and stable LKCD examples as the objects of primary interest.
Section~\ref{sec:ggc} discusses the ID class known as generalised gamma convolutions (GGCs),
followed by the stable GGC LKCD in Section~\ref{sec:stableGGCLKCD}.
With the preparatory background in place, Section~\ref{sec:convolution} moves to the convolution of two gamma densities
with several associated LKCD examples.
This is the first part of the core message of the paper.
The second part in  Section~\ref{sec:stablemixtures} introduces mixtures of stable densities, 
with both stable and gamma mixing densities as  examples.
The former allows the generation of new stable densities from given instances. 
The latter gives the fractional gamma density, which is discussed at length.
We discover  a novel integral  representation of the fractional gamma density for $\alpha=1/2$ in terms of
the parabolic cylinder function.
All examples are presented as LKCDs.
This is followed by a discussion in Section~\ref{sec:discussion}, and concluding  remarks and pointers to future work in
Section~\ref{sec:conclusion}.


\subsection{Stable distribution}
 \label{sec:stable} 
The $\alpha$-stable distribution  for $0<\alpha<1$, defined  on the positive  half-line,  has density $f_\alpha(x), \, x\ge 0$
with  Laplace transform 
\begin{align}
\widetilde f_\alpha(s) &= \exp(-s^\alpha) 
\label{eq:laplacestable}
\end{align}
The $\alpha$-stable distribution is of interest in probability theory and various applications. 
In physics, $\widetilde f_\alpha(s)$ is often referred to  as the stretched exponential or the Kohlrausch function
(Berberan-Santos~{\it et al.}~\cite{BERBERANSANTOS2005171}, Penson and G\'{o}rska~\cite{PensonGorska}).
It is intimately associated with  relaxation and diffusion phenomena.
To paraphrase~\cite{PensonGorska}, $f_\alpha(x)$ arises in 
condensed and soft matter physics, geophysics, 
meteorology, economics, fractional kinetics:
\textquote{For instance, the value $\alpha= 1/4$ is thought to describe mechanical and dielectric properties of glassy polymers. 
It is also confirmed that the same value of $\alpha$  is relevant for a statistical description of subrecoil laser cooling}.

Yet the functional form of $f_\alpha(x)$ is  elusive.
Pollard~\cite{Pollard} showed that Laplace inversion 
gives the rather forbidding infinite series 
\begin{align}
f_\alpha(x) &= -\frac{1}{\pi}\sum_{k=0}^\infty \frac{(-1)^k}{k!} \sin(\pi k \alpha)\frac{\Gamma(k \alpha+1)}{x^{k \alpha +1}}
\label{eq:pollardsum}
\end{align}
Feller~\cite{Feller2} (p581) derived the expression using the Fourier transform, 
along with another expression  for $1<\alpha<2$.
(We shall not discuss here  stable distributions on the real line indexed by $1\le \alpha\le2$, with known closed form  only for
the Cauchy distribution for $\alpha=1$ and the Gaussian distribution for $\alpha=2$.)

For $\alpha=1/2$, the series representation~(\ref{eq:pollardsum}) reduces to the simple  form
\begin{align}
f_{\tfrac{1}{2}}(x) &= \frac{1}{2\sqrt{\pi}}\,x^{-3/2}\,e^{-1/4x}
\label{eq:alfahalf}
\end{align}
Other  forms  can be inferred from~(\ref{eq:pollardsum}), 
such as $f_{1/3}$  in terms of $K_{1/3}$, the modified Bessel function of the second kind of order 1/3.
Forms for  rational  $\alpha$ are typically cast in terms  of hypergeometric functions or the  allied  Whittaker functions,
{\it e.g.}\ $\alpha=\{2/3, 1/4, 3/4\}$~(Barkai~\cite{Barkai},  Penson and G\'{o}rska~\cite{PensonGorska}, 
Scher and Montroll~\cite{ScherMontroll}).
Indeed, ~\cite{PensonGorska,ScherMontroll} state that, for any rational $\alpha=l/k$ ($0<l<k$), 
$f_\alpha$ may be expressed as a finite sum of generalised hypergeometric functions.
But generality often comes at the cost of simplicity.
Hypergeometric functions are flexible series representations that are not routinely encountered mathematical objects, 
even though they yield many common functions as particular cases.

Amongst other things, we will discuss a simple and  known integral representation  of $f_{\alpha\beta}$ in terms  of 
$f_\alpha$ and  $f_\beta$.
In particular, we shall  infer $f_{1/4}$ from $f_{1/2}$
as an integral representation  instead of the hypergeometric representation of $f_{1/4}$ given in~\cite{Barkai, PensonGorska}.

In Section~\ref{sec:LKCD} we shall  motivate the stable distribution as an intimate relative of the gamma distribution 
rather than as a standalone object.
To that end and beyond, we discuss  next the concept of infinitely divisible (ID) distributions.

\section{Infinitely Divisible Distributions}
\label{sec:ID}
The theory of infinitely divisible distributions summarised here is well-known and can be found in several probabilistic
texts such as Feller~\cite{Feller2}, Kingman~\cite{KingmanBook}, Steutel and van Harn~\cite{SteutelvanHarn}.
Our contribution  is  a  commutative diagram representation that, in our view, offers a helpful  visual summary of the  
theoretical framework.


By way of  basic motivation, consider a set of points that are  randomly scattered over some domain. 
In practical  application, the domain might be  an interval in time or a region in space.
A point might be an event in time like  a vehicle crossing a bridge in sparse traffic, an isolated day of rain 
or, in a spatial  context, a  point source at some location in the sky.
 In addition, each point carries a random positive additive attribute, 
such as the mass of the vehicle, the amount of rainfall on the given day or  the brightness of the point source.
In each case, we may meaningfully speak of the  total vehicle mass  that the bridge bears in a day, the   rainfall in a month
or the  brightness of the  patch of sky by simply adding up the respective attributes over all  point occurrences within the 
specified domain. 

More abstractly, 
let $n$, the number of point occurrences in a specified domain,  be governed by a 
Poisson distribution 
with mean rate $\mu$ (typically the size of the domain).
Let each point $i$ have an associated attribute  $X_i$ where the $\{X_i: i=1\ldots n\}$ are  
independent, identically distributed  positive  (nonnegative) variables governed by a common distribution
with density  $\ell(x)$.
Then, as is well-known, the sum  $X=X_1+X_2+\dots+X_n$ is governed by the  density $\Pr(x|n)= \ell^{n\star}(x)$,
where $\ell^{n\star}$ is the $n$-fold Laplace convolution of $\ell$ (where $\ell^{1\star} \equiv \ell$).
 For $n=0$, $X\equiv 0$,  so that $\Pr(x|n=0)\equiv \ell^{0*}(x)=\delta(x)$ is an atom at $x=0$.

The joint distribution of $X$ and $n$ is 
\begin{align}
\Pr(x,n|\mu) = \Pr(x|n)\Pr(n|\mu) &= \ell^{n\star}(x)\, e^{-\mu} \,\frac{\mu^n}{n!}
\label{eq:prxn}
\intertext{Hence the unconditional distribution of $X$ is}
\Pr(x|\mu) = \sum_{n=0}^\infty \Pr(x|n)\Pr(n|\mu) &=  e^{-\mu}\sum_{n=0}^\infty  \frac{\mu^n}{n!}\, \ell^{n\star}(x)
\label{eq:prx}
\end{align}
We shall  also write this as  $f(x|\mu)$, which is 
the density of what is known as the compound Poisson distribution
that we shall  denote by $\cC\cP(\mu,\ell)$.
The Laplace  transform  of $f(x|\mu)$ is
\begin{align}
{\scrL}\{f\}(s) \equiv {\widetilde f}(s|\mu) &= \int_0^\infty e^{-sx}\ f(x|\mu) dx
\label{eq:lpcf}
\end{align}
Similarly,  $\widetilde \ell(s)$ is the  Laplace  transform of $\ell(x)$. 
Since, by the convolution theorem, the Laplace transform of a convolution of functions  is a product of their respective  Laplace transforms,
$\scrL\{\ell^{n\star}\}(s)=\widetilde\ell\,^{n}(s)$.
Hence the Laplace transform of~(\ref{eq:prx}) is
\begin{align}
\widetilde f(s|\mu) =  e^{-\mu}\sum_{n=0}^\infty  \frac{\mu^n}{n!}\, \widetilde\ell\,^{n}(s)
&= \exp\{-\mu(1-\widetilde\ell(s))\} \\
&= \exp\left\{- \mu\int_0^\infty \left(1-e^{-sx}\right)\ell(x)dx\right\}
\label{eq:LK}
\end{align}
where~(\ref{eq:LK}) follows from $\ell(x)$ being a  density that is normalised
(at least at this stage of the discussion). 
A distribution with Laplace transform~(\ref{eq:LK}) is said to be infinitely divisible because any $n^{\rm th}$ root of~(\ref{eq:LK}) 
is the same expression with $\mu$ replaced by $\mu/n$,
{\it i.e.}\ the $n^{\rm th}$ root is the Laplace transform of the probability distribution with density $f(x|\tfrac{\mu}{n})$.
The form (\ref{eq:LK})  is the  celebrated  L\'{e}vy-Khintchine representation of an infinitely divisible distribution on positive additive variables,
with the definition of $\ell(x)$, known as the L\'{e}vy density,  broadened beyond a normalised density. 
We may write  the  Laplace exponent as
\begin{align}
  \psi(s)\equiv\int_0^\infty \left(1-e^{-sx}\right)\ell(x)dx 
  &=  \int_0^\infty \int_0^s e^{-xt}dt\, x\ell(x)dx \\
  &=  \int_0^s  \int_0^\infty e^{-xt} \rho(x)dx\, dt \\
  &=  \int_0^s \widetilde\rho(t)dt
 \label{eq:LKexp}
\end{align}
where $\widetilde\rho(s)$ is the Laplace transform of $\rho(x)\equiv x\ell(x)$.
In light of~(\ref{eq:LKexp}),~(\ref{eq:LK}) becomes 
\begin{align}
\widetilde f(s|\mu) &= \exp\left(-\mu\,\psi(s)\right)
\label{eq:LKa}
\end{align}

Hence it  is  $\widetilde\rho(s)$ 
that actually needs to exist rather  than $\widetilde\ell(s)$.
 The compound Poisson representation $\cC\cP(\mu,\ell)$ need not strictly exist for infinite  divisibility to hold.
Differentiating~(\ref{eq:LKa}) gives 
\begin{align}
\widetilde f\,^\prime(s|\mu) &= -\mu\widetilde\rho(s)\widetilde f(s|\mu) \\
\implies\qquad 
\mu\,\widetilde\rho(s) &= - \frac{\widetilde f\,^\prime(s|\mu)}{\widetilde f(s|\mu)}
\label{eq:derivLK} 
\end{align}
which is invariant under scaling $f(x|\mu)\to Cf(x|\mu)$ for any constant $C>0$.
An  equivalent expression arises from the limiting process 
\begin{align}
\lim_{n\to\infty}n \widetilde f\,^\prime(s|\tfrac{\mu}{n})  &=  -\lim_{n\to\infty} \mu\widetilde\rho(s)\widetilde f(s|\tfrac{\mu}{n})  
     = -\mu \widetilde\rho(s) \widetilde f(s|0)
\label{eq:limLK} \\
\implies\quad 
\mu\,\widetilde\rho(s) &= -\lim_{n\to\infty}n\,\widetilde f\,^\prime(s|\tfrac{\mu}{n})/\widetilde f(s|0)
    = -\lim_{n\to\infty}n\,\widetilde f\,^\prime(s|\tfrac{\mu}{n})
\label{eq:limderivLK} 
\end{align}
given that, by~(\ref{eq:LKa}), $\lim_{n\to\infty} \widetilde f(s|\tfrac{\mu}{n})= \widetilde f(s|0)=1$.
Invariance  under $f(x|\mu)\to Cf(x|\mu)$ is preserved because, correspondingly, $\widetilde f(s|0)=1\to \widetilde f(s|0)=C$. 
Since $-\widetilde f\,^\prime(s|\mu)$ is the Laplace  transform of $xf(x|\mu)$, it follows that
\begin{align}
\mu \rho(x) &= \lim_{n\to\infty}n xf(x|\tfrac{\mu}{n})/\widetilde f(s|0) = \lim_{n\to\infty}n xf(x|\tfrac{\mu}{n})
\quad (\widetilde f(s|0)=1)
\label{eq:limnxf} 
\end{align}

An alternative approach to the foregoing starts from Bernstein's theorem~\cite{Feller2}~(p439), 
which  states that a function $f(x)$ is a density
if and only if its  Laplace transform $\widetilde f(s)$ is completely monotone, {\it i.e.}\ $(-1)^n\widetilde f^{(n)}(s)\ge0$, $n\ge0$ 
($f(x)$ is a probability density if, in addition, $\widetilde f(0)=1$).
Infinite divisibility  of a density $f(x)$ is the case where $-\widetilde f\,^\prime(s)/\widetilde f(s)$ is also the  Laplace transform of a density.
Hence an alternative definition of infinite divisibility  is that $f(x)$ is the  density of an
infinitely divisible distribution  if and only if  both $\widetilde f(s)$ and $-\widetilde f\,^\prime(s)/\widetilde f(s)$ are  completely monotone.
This more abstract definition  is consistent with the compound Poisson approach although it does not directly assume it.
In the compound Poisson construction, 
$\ell(x)$ and therefore $\rho(x)=x\ell(x)$ is an assigned density from the outset so that $\widetilde\rho(s)$ is necessarily completely monotone.

The Laplace convolution of two or more densities will arise repeatedly in the discussion that follows.
Invoking the convolution theorem once more, it is straightforward to see that the convolution of  ID densities
is also an ID density whose L\'{e}vy density is the sum of L\'{e}vy densities of the  convolution components.

\section{L\'{e}vy-Khintchine Commutative Diagram}
\label{sec:LKCD}
\begin{figure}[tbh]
\begin{center}
\begin{tikzpicture}[auto,scale=2.7]
\newcommand*{\size}{\scriptsize}%
\newcommand*{\gap}{.2ex}%
\newcommand*{\width}{3.5}%
\newcommand*{\height}{2.0}%

\node (P) at (0,0)  {$f(x|\mu)$};
\node (Q) at ($(P)+(\width,0)$)  {$\widetilde f(s|\mu)$};
\node (B) at ($(P)-(0,\height)$) {$\mu\,\rho(x)$};
\node (C) at  ($(B)+(\width,0)$) {$\mu\,\widetilde\rho(s)$};    
\draw[Myarrow] ([yshift =  \gap]P.east)  --  node[above] {\size $\scrL$} ([yshift = \gap]Q.west) ;
\draw[Myarrow] ([yshift = -\gap]Q.west) -- node[below] {\size ${\scrL}^{-1}$} ([yshift = -\gap]P.east); 
\draw[Myarrow]([xshift =  \gap]Q.south) -- 
     node[pos=0.35] {\size $-\dfrac{\widetilde f^\prime(s|\mu)}{\widetilde f(s|\mu)}$ \textrm{ or}}
     node[pos=0.65] {\size $-\displaystyle\lim_{n\to\infty}n\widetilde f\,^\prime(s|\tfrac{\mu}{n})$}
([xshift =  \gap]C.north);
\draw[Myarrow]([xshift = -\gap]C.north) -- node {\size $\exp{\left\{-\mu\psi(s)\right\}}$}([xshift = -\gap]Q.south);
\draw[Myarrow,dashed] ([xshift = -\gap]B.north) -- node {\size ${{\cal CP}(\mu,\ell)}$} ([xshift =  -\gap]P.south);
\draw[Myarrow] ([xshift =  \gap]P.south) -- node {\size ${{\displaystyle\lim_{n\to\infty}}n xf(x|\frac{\mu}{n})}$} ([xshift = \gap]B.north);
\draw[Myarrow] ([yshift =  \gap]B.east)  --  node[above] {\size $\scrL$} ([yshift =  \gap]C.west) ;
\draw[Myarrow] ([yshift = -\gap]C.west) --  node[below] {\size ${\scrL}^{-1}$}  ([yshift = -\gap]B.east); 
\draw[Myarrow] (B.north east)  --   
                node[sloped, anchor=east, centered, fill=white] {\size $\exp\left\{- \mu\displaystyle\int_0^\infty (1-e^{-sx})\,\ell(x)dx\right\}$}
                (Q.south west);
\end{tikzpicture}
\end{center}
\caption{L\'{e}vy-Khintchine Commutative Diagram (LKCD). 
$\scrL$ is the Laplace transform and ${\cC\cP(\mu,\ell)}$ ($x\ell(x)=\rho(x)$) is the compound Poisson construction 
(dotted because it may formally be undefined).
$\psi(s)$ is the (definite or indefinite) integral of $\widetilde\rho(s)$. 
The direct  L\'{e}vy-Khintchine relation  is the  diagonal  from bottom left to top right.
It is  equivalent to a composition of  transitions along the axes: 
\textquote{east then north} or (if ${\cC\cP(\mu,\ell)}$ exists) \textquote{north then east}.}
\label{fig:LKCD}
\end{figure}
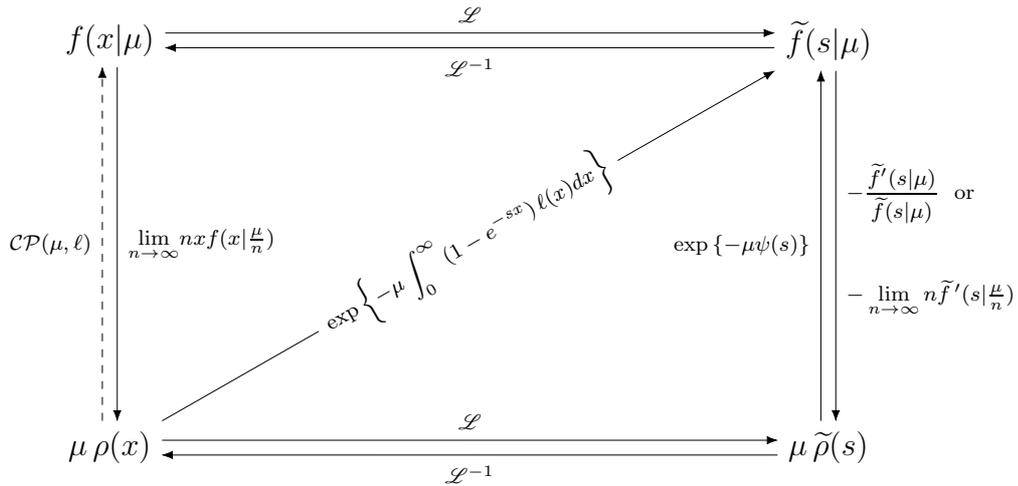

We summarise the objects and relationships amongst them in a   graphic that we refer to as the 
 L\'{e}vy-Khintchine commutative diagram (LKCD), shown in Figure~\ref{fig:LKCD}.
 The compound Poisson relation is  dotted to accommodate  the observation above that it may formally be undefined 
 despite the existence of all nodes of the LKCD
 (this  may differ from the conventional interpretation of a dotted arrow in category theory).
  
The commutative diagram illustrates, at a glance, the concept of infinite divisibility, the objects involved and the relationships amongst them.
 Having assigned or constructed one of  the four possible nodes, we may then seek to populate the other nodes by 
 following a path of  relationships best suited to the task.
Although the LKCD  is  primarily an organising principle rather than a  discovery tool ({\it i.e.}\ a mechanism to construct new ID distributions),  
in our view the visual representation facilitates both  the description and construction of  ID distributions.

We shall often speak of the upper level of the LKCD as the ID density level and
the lower level as the  L\'{e}vy density level.
As noted earlier, any multiplicative constant at the upper level is `forgotten' upon descent to the lower level.
The density $\rho(x)$ in the lower level may itself be ID,  
in which case the LKCD  can be extended downward to form a two-level `ladder' where the bottom rung is the
L\'{e}vy density level of the middle layer, which is in turn the L\'{e}vy density level  of the top layer.
Furthermore, since $f(x|\mu)$ is a density, of necessity, the layer $(f(x|\mu), \widetilde f(s|\mu))$ can   be treated as
a L\'{e}vy density level for a higher level ID density layer, thereby extending the LKCD ladder upward by another rung.
In principle, such upward growth  of the LKCD can be repeated indefinitely, 
although analytic expressions for the  ID densities thus generated  may become increasingly elusive.

We give LKCD  examples for two densities that are central to the rest of our discussion.

\subsection{Gamma  LKCD}
\label{sec:gammaLKCD}
The gamma density $g_{\mu,\lambda}(x)\equiv f_\lambda(x|\mu)$ and its LKCD are, respectively~(\ref{eq:gamma}) 
and~(\ref{eq:gammaLKCD})
\begin{align}
g_{\mu,\lambda}(x) &= \frac{\lambda^\mu}{\Gamma(\mu)}\,x^{\mu-1}e^{-\lambda x} 
    = \frac{\mu\, \lambda^\mu}{\Gamma(1+\mu)}\,x^{\mu-1}e^{-\lambda x}\qquad \mu,\, \lambda>0
\label{eq:gamma}
\end{align}
\begin{equation}
\begin{tikzpicture}[auto,scale=1.75, baseline=(current  bounding  box.center)]
\newcommand*{\size}{\scriptsize}%
\newcommand*{\gap}{.2ex}%
\newcommand*{\width}{2.25}%
\newcommand*{\height}{1.25}%

\node (P) at (0,0)  {$\dfrac{\lambda^\mu}{\Gamma(\mu)}\,x^{\mu-1}e^{-\lambda x}$};
\node (Q) at ($(P)+(\width,0)$) {$\dfrac{\lambda^\mu}{(\lambda+s)^\mu}$};
\node (B) at ($(P)-(0,\height)$) {$\mu\, e^{-\lambda x}$};
\node (C) at  ($(B)+(\width,0)$) {$\dfrac{\mu}{\lambda+s}$};   

\draw[Myarrow] ([yshift =  \gap]P.east)   --   ([yshift =  \gap]Q.west) ;
\draw[Myarrow] ([yshift = -\gap]Q.west)  --   ([yshift = -\gap]P.east); 
\draw[Myarrow]([xshift  =  \gap]Q.south) --  ([xshift =  \gap]C.north);
\draw[Myarrow]([xshift  = -\gap]C.north) --   ([xshift = -\gap]Q.south);
\draw[Myarrow] ([xshift =  \gap]P.south) -- ([xshift =  \gap]B.north);
\draw[Myarrow] ([yshift =  \gap]B.east)  --  ([yshift =  \gap]C.west) ;
\draw[Myarrow] ([yshift = -\gap]C.west) --   ([yshift = -\gap]B.east); 

\node (P) at (\width*1.75,0)  {$\dfrac{1}{\Gamma(\mu)}x^{\mu-1}$};
\node (Q) at ($(P)+(\width,0)$) {$\dfrac{1}{s^\mu}$};
\node (B) at ($(P)-(0,\height)$) {$\mu$};
\node (C) at  ($(B)+(\width,0)$) {$\dfrac{\mu}{s}$};   
\node (A) at ($(P)+(\width/2,-\height/2)$) {$\lambda=0$}; 

\draw[Myarrow] ([yshift =  \gap]P.east)   --   ([yshift =  \gap]Q.west) ;
\draw[Myarrow] ([yshift = -\gap]Q.west)  --   ([yshift = -\gap]P.east); 
\draw[Myarrow]([xshift  =  \gap]Q.south) --  ([xshift =  \gap]C.north);
\draw[Myarrow]([xshift  = -\gap]C.north) --   ([xshift = -\gap]Q.south);
\draw[Myarrow] ([xshift =  \gap]P.south) -- ([xshift =  \gap]B.north);
\draw[Myarrow] ([yshift =  \gap]B.east)  --  ([yshift =  \gap]C.west) ;
\draw[Myarrow] ([yshift = -\gap]C.west) --   ([yshift = -\gap]B.east); 
\end{tikzpicture}
\label{eq:gammaLKCD}
\end{equation}

The compound Poisson construction  is not defined 
but all other mappings are well-defined.
In particular, the second form of~(\ref{eq:gamma}) makes it clear that the limit~(\ref{eq:limnxf}) is
\begin{align}
\lim_{n\to\infty}n\, x\, f_\lambda(x|\tfrac{\mu}{n}) &= \mu \rho_\lambda(x) = \mu\, e^{-\lambda x}
\label{eq:limgamma}
\end{align}
We  may omit the  normalising factor $\lambda^\mu$, 
thereby allowing the case $\lambda=0$ to be well-defined, 
as shown in  the second frame of~(\ref{eq:gammaLKCD}).
The density is no longer finite ($s^{-\mu}$ is not  defined at $s=0$), but the LKCD representation remains valid. 
Henceforth we shall  routinely omit $\lambda^\mu$ in the definition of the density $g_{\mu,\lambda}(x)$,
thereby making it valid for $\mu>0$ and $\lambda\ge 0$.

\subsection{Stable  LKCD}
\label{sec:stableLKCD}
We now consider upward extension  of the gamma LKCD~(\ref{eq:gammaLKCD})
 for  $\mu=1-\alpha$ ( $0<\alpha<1$) and $\lambda=0$.
We  reuse $\mu>0$ as a multiplicative  parameter that  will be a scale parameter for the layer above.
We further multiply by $\alpha$ to generate $\rho_\alpha(x) = \alpha\,x^{-\alpha}/\Gamma(1-\alpha)= \alpha \, g_{1-\alpha,0}(x)$.
As a gamma density, $\rho_\alpha(x)$ is ID, exactly as discussed above and shown in the unshaded LKCD below.
But 
$\rho_\alpha(x)/x$ is  itself  the L\'{e}vy density
 of a higher ID density $f_\alpha(x)$, as shown in the upper shaded extension of the  LKCD.
\begin{equation}
\begin{tikzpicture}[auto,scale=1.75, baseline=(current  bounding  box.center)]
\newcommand*{\size}{\scriptsize}%
\newcommand*{\gap}{.2ex}%
\newcommand*{\width}{2.0}%
\newcommand*{\height}{1.25}%

\node (P) at (0,0)  {$f_{\alpha}(x|\mu)$};
\node (Q) at ($(P)+(\width,0)$)   {$\exp(-\mu s^\alpha)$};
\node (B) at ($(P)-(0,\height)$) {$\mu\alpha\,\dfrac{\,x^{-\alpha}}{\Gamma(1-\alpha)}$};
\node (C) at ($(B)+(\width,0)$) {$\mu\alpha\, s^{\alpha-1}$};    
\node (X) at ($(B)-(0,\height)$) {$1-\alpha$};
\node (Y) at ($(X)+(\width,0)$) {$\dfrac{1-\alpha}{s}$};    

\node (B1) at ($(B)+(0,7.5*\gap)$) {}; 
\node (C1) at ($(B1)+(\width,0)$) {}; 

\draw[Myarrow] ([yshift =  \gap]P.east)   --   ([yshift =  \gap]Q.west) ;
\draw[Myarrow] ([yshift = -\gap]Q.west)  --   ([yshift = -\gap]P.east); 
\draw[Myarrow]([xshift  =  \gap]Q.south) --  ([xshift =  \gap]C1.north);
\draw[Myarrow]([xshift  = -\gap]C1.north) --   ([xshift = -\gap]Q.south);
\draw[Myarrow,dashed] ([xshift = -\gap]B1.north) -- ([xshift =  -\gap]P.south);
\draw[Myarrow] ([xshift =  \gap]P.south) -- ([xshift =  \gap]B1.north);
\draw[Myarrow] ([yshift =  \gap]B.east)  --  ([yshift =  \gap]C.west) ;
\draw[Myarrow] ([yshift = -\gap]C.west) --   ([yshift = -\gap]B.east); 
\draw[Myarrow] ([yshift =  \gap]X.east)  --  ([yshift =  \gap]Y.west) ;
\draw[Myarrow] ([yshift = -\gap]Y.west) --   ([yshift = -\gap]X.east); 
\draw[Myarrow,dashed] ([xshift = -\gap]X.north) -- ([xshift =  -\gap]B.south);
\draw[Myarrow] ([xshift =  \gap]B.south) -- ([xshift =  \gap]X.north);
\draw[Myarrow]([xshift  = -\gap]Y.north) --   ([xshift = -\gap]C.south);
\draw[Myarrow] ([xshift =  \gap]C.south) -- ([xshift =  \gap]Y.north);

\fill[opacity=0.15]  ($(P)+(-\width*0.3,\height*0.2)$) rectangle ($(C)+(\width*0.3,+\height*0.22)$);

\node (P) at (\width*2.0,0)  {$\dfrac{\mu}{2\sqrt{\pi x^3}}\,e^{-\mu^2/4x}$};
\node (Q) at ($(P)+(\width,0)$) {$\exp(-\mu \sqrt{s})$};
\node (B) at ($(P)-(0,\height)$) {$\dfrac{\mu}{2\sqrt{\pi x}}$};
\node (C) at  ($(B)+(\width,0)$) {$\dfrac{\mu}{2\sqrt{s}}$};   
\node (X) at ($(B)-(0,\height)$) {$\dfrac{1}{2}$};
\node (Y) at ($(X)+(\width,0)$) {$\dfrac{1}{2s}$};  
\node (A) at ($(P)+(\width/2,-\height*1.5)$) {$\alpha=\tfrac{1}{2}$}; 

\node (B1) at ($(B)+(0,7.5*\gap)$) {}; 
\node (C1) at ($(B1)+(\width,0)$) {}; 

\draw[Myarrow] ([yshift =  \gap]P.east)   --   ([yshift =  \gap]Q.west) ;
\draw[Myarrow] ([yshift = -\gap]Q.west)  --   ([yshift = -\gap]P.east); 
\draw[Myarrow]([xshift  =  \gap]Q.south) --  ([xshift =  \gap]C1.north);
\draw[Myarrow]([xshift  = -\gap]C1.north) --   ([xshift = -\gap]Q.south);
\draw[Myarrow,dashed] ([xshift = -\gap]B1.north) -- ([xshift =  -\gap]P.south);
\draw[Myarrow] ([xshift =  \gap]P.south) -- ([xshift =  \gap]B1.north);
\draw[Myarrow] ([yshift =  \gap]B.east)  --  ([yshift =  \gap]C.west) ;
\draw[Myarrow] ([yshift = -\gap]C.west) --   ([yshift = -\gap]B.east); 
\draw[Myarrow] ([yshift =  \gap]X.east)  --  ([yshift =  \gap]Y.west) ;
\draw[Myarrow] ([yshift = -\gap]Y.west) --   ([yshift = -\gap]X.east); 
\draw[Myarrow,dashed] ([xshift = -\gap]X.north) -- ([xshift =  -\gap]B.south);
\draw[Myarrow] ([xshift =  \gap]B.south) -- ([xshift =  \gap]X.north);
\draw[Myarrow]([xshift  = -\gap]Y.north) --   ([xshift = -\gap]C.south);
\draw[Myarrow] ([xshift =  \gap]C.south) -- ([xshift =  \gap]Y.north);

\fill[opacity=0.15]  ($(P)+(-\width*0.4,\height*0.2)$) rectangle ($(C)+(\width*0.3,+\height*0.22)$);

\end{tikzpicture}
\label{eq:3levelstableLKCD}
\end{equation}
The shaded top layer, induced by the middle layer as its  L\'{e}vy density level, is the stable density $f_\alpha(x|\mu)$ 
with Laplace  transform $\widetilde f_\alpha(s|\mu)=\exp(-\mu s^\alpha)$.
All nodes are filled from  knowledge of the Laplace transform except for  $f_\alpha(x|\mu)$  itself, 
for which there is no known simple simple closed form expression for general $\alpha$.
The LKCD for $\alpha=1/2$, for which $f_{1/2}(x|\mu)$ is known, 
is shown on the right of~(\ref{eq:3levelstableLKCD}).

The Stable LKCD can be further extended, rather naturally, to the left,
to form what is known as a generalised gamma convolution.

\section{Generalised Gamma Convolution}
\label{sec:ggc}
In keeping with the  initial definition, the minimal LKCD for an ID density $f(x|\mu)$ is the 
rectangle on the right of~(\ref{eq:GGCLKCD}).
\begin{equation}
\begin{tikzpicture}[auto,scale=2.0, baseline=(current  bounding  box.center)]
\newcommand*{\size}{\scriptsize}%
\newcommand*{\gap}{.2ex}%
\newcommand*{\width}{1.75}%
\newcommand*{\height}{1.0}%

\node (O) at (0,0) {}; 
\node (P) at ($(O)+(\width,0)$) {$f(x|\mu)$};
\node (Q) at ($(P)+(\width,0)$) {$\widetilde f(s|\mu)$};
\node (A) at ($(O)-(0,\height)$) {$\mu\, u(t)$};
\node (B) at ($(A)+(\width,0)$) {$\mu\,\rho(x)$};
\node (C) at ($(B)+(\width,0)$) {$\mu\,\widetilde\rho(s)$};    

\draw[Myarrow] ([yshift =  \gap]P.east)   --   ([yshift =  \gap]Q.west) ;
\draw[Myarrow] ([yshift = -\gap]Q.west)  --   ([yshift = -\gap]P.east); 
\draw[Myarrow]([xshift  =  \gap]Q.south) --  ([xshift =  \gap]C.north);
\draw[Myarrow]([xshift  = -\gap]C.north) --   ([xshift = -\gap]Q.south);
\draw[Myarrow,dashed] ([xshift = -\gap]B.north) -- ([xshift =  -\gap]P.south);
\draw[Myarrow] ([xshift =  \gap]P.south) -- ([xshift =  \gap]B.north);
\draw[Myarrow] ([yshift =  \gap]A.east)  --  ([yshift =  \gap]B.west) ;
\draw[Myarrow] ([yshift = -\gap]B.west) --   ([yshift = -\gap]A.east); 
\draw[Myarrow] ([yshift =  \gap]B.east)  --  ([yshift =  \gap]C.west) ;
\draw[Myarrow] ([yshift = -\gap]C.west) --   ([yshift = -\gap]B.east); 

\end{tikzpicture}
\label{eq:GGCLKCD}
\end{equation}
If, in addition, there exists a density $u(t)$  with Laplace transform $\rho(x)$
(equivalently, if $\rho(x)$ is completely monotone) then $f(x|\mu)$ is said to be a 
generalised gamma convolution (GGC).
We refer to the extended LKCD~(\ref{eq:GGCLKCD}), including $u(t)$, as a GGC LKCD.
The GGC name derives from  the fact that an arbitrary  sum of delta functions $u(t)=\sum_i u_i\, \delta(t-\lambda_i)$ has Laplace transform 
$\widetilde u(x) \equiv\rho(x)= \sum_i u_i\, \exp{-\lambda_i x}$, whose $f(x|\mu)$ is the
convolution of as many gamma densities corresponding to each term in the $\rho(x)$ sum of exponentials.
In general, $\rho(x)$ need not be a sum of exponentials, any density $\rho(x)$ that is itself the Laplace  transform of another density $u(t)$
generates an $f(x|\mu)$ known as a GGC.
For example, $\rho(x)$ might itself be a gamma density, as is the case for the stable density, which is thus an instance of a GGC.

The GGC class  was introduced by Thorin~\cite{Thorin77a, Thorin77b} 
as he sought to prove the infinite divisibility of the Pareto  and the log-normal distributions. 
The general density $u(t)$, if it exists,  is known as the Thorin density.
GGCs were  subsequently studied in depth in the book by Bondesson~\cite{Bondesson}. 


An alternative but equivalent motivation, as given  in the survey of GGCs by James~{\it et al.}~\cite{JamesRoynetteYor}, 
proceeds as follows.
First introduce the concept of self-decomposability, which may be compactly defined as follows (Bondesson~\cite{Bondesson},~p18):
\begin{quote}
An ID distribution $f(x)$ on $\mathbb{R}_+$ is self-decomposable if and only if it has a  L\'{e}vy density $\ell(x)$  such that
$\rho(x)=x\ell(x)$ is decreasing.
\end{quote}
If, in addition to $\rho(x)$ decreasing, there also exists a density $u(t)$ on  $\mathbb{R}_+$ such that $\rho(x)$ is the Laplace transform 
of $u(t)$, then  $f(x)$ is GGC with $u(t)$ as the Thorin density
(equivalently, $\rho(x)$ has to be  decreasing {\it and} completely monotone for $f(x)$ to be GGC).
The following class hierarchy holds:
\vspace{-10pt}
\begin{myenumerate}
\item Let $\cI$ be the class of ID distributions on $\mathbb{R}_+$ 
\item Let $\cS$ be  the class of self-decomposable distributions on $\mathbb{R}_+$
\item Let $\cG$ be  the class of generalised gamma convolutions (necessarily defined on $\mathbb{R}_+$)
\end{myenumerate}
\vspace{-5pt}
Then $\cG\subset\cS\subset\cI$, as illustrated in Bondesson~\cite{Bondesson}~p4 (with additional subclasses of $\cG$ that
we shall not explicitly discuss  here, such as the hyperbolically  completely monotone distributions).
The gamma distribution is, of course, GGC and so is the closely allied stable distribution, as discussed  in the next section.
Outside GGC, we shall not explore any other particular cases of  ID distributions on $\mathbb{R}_+$, such as, say, 
distributions that may be self-decomposable but not GGC.

There is a fundamental relationship between the $\widetilde\rho(s)$ and $u(t)$.
Given $\widetilde\rho(s)$, Bondesson~\cite{Bondesson} (p33, Inversion Theorem) used contour integration to show that
 \begin{align}
\mu\,u(t) = \frac{1}{\pi}\,{\rm Im}\, \mu\,\widetilde\rho(-t) 
                 &= \frac{1}{\pi}\,{\rm Im}\,\frac{ \widetilde f\,^\prime(-t|\mu)}{\widetilde f(-t|\mu)}                  
\label{eq:derivGGCinversion} \
\intertext{which, by the limiting rule~(\ref{eq:limderivLK}) above, is equivalent to}
\mu\,u(t) = \frac{1}{\pi}\,{\rm Im}\, \mu\,\widetilde\rho(-t) 
                &=  \frac{1}{\pi}\,{\rm Im}\, \lim_{n\to\infty}n\widetilde f\,^\prime(-t|\tfrac{\mu}{n})
\label{eq:limGGCinversion} 
 \end{align}
 Bondesson's proof relied on the theory of Pick functions, defined as functions that are analytic on the upper complex half-plane.
 We shall give a simple  derivation ({\it i.e.}\ without invoking contour integration) of the Inversion Theorem for the stable distribution. 
 We shall use Pollard's infinite series representation of 
 $f_\alpha(x|\mu)$  as the basis for discussion.  
 
 
\section{The Stable GGC Commutative Diagram }
\label{sec:stableGGCLKCD}

Pollard~\cite{Pollard} used contour integration to derive  the following integral representation
\begin{align}
f_\alpha(x|\mu=1) \equiv f_\alpha(x)  &= \frac{1}{\pi} {\rm Im} \int_0^\infty e^{-xt}\,e^{-(te^{-i\pi})^\alpha}\,dt 
\label{eq:pollardint1} \\
  &= \frac{1}{\pi}\int_0^\infty e^{-xt}\,e^{-t^\alpha \cos\pi\alpha} \sin(t^\alpha \sin\pi\alpha)\, dt
\label{eq:pollardint2}
\end{align}
The change of variable $x\to\mu^{-1/\alpha}x$ takes $f_\alpha(x)dx$ to 
$f_\alpha(x|\mu)dx\equiv f_\alpha(\mu^{-1/\alpha}x)\mu^{-1/\alpha}dx$ 
and the   Laplace transform $\widetilde f_\alpha(s)=\exp(-s^\alpha)$ to  $\widetilde f_\alpha(\mu^{1/\alpha}s)= \exp(-\mu s^\alpha)$.
Thus, starting with  the infinite series~(\ref{eq:pollardsum}), we readily arrive at  the equivalent   integral representation for 
$f_\alpha(x|\mu)$: 
\begin{align}
f_\alpha(x|\mu) &= -\frac{1}{\pi}\sum_{k=0}^\infty \frac{(-1)^k\mu^k}{k!} \sin(\pi k \alpha)\frac{\Gamma(k \alpha+1)}{x^{k \alpha +1}} 
\label{eq:pollardsum2} \\
    &= \frac{1}{\pi}\sum_{k=0}^\infty \frac{(-1)^k\mu^k}{k!} \sin(-\pi k \alpha)\int_0^\infty e^{-xt} t^{k\alpha} dt \\
    &= \frac{1}{\pi}\,{\rm Im} \int_0^\infty e^{-xt} \sum_{k=0}^\infty \frac{(-1)^k}{k!} (\mu\, e^{-i\pi\alpha} t^\alpha)^k\,dt \\
    &= \frac{1}{\pi}\,{\rm Im} \int_0^\infty e^{-xt} e^{-\mu(e^{-i\pi}t)^\alpha} \, dt \\
    &= \frac{1}{\pi}\,{\rm Im} \int_0^\infty e^{-xt}\, \widetilde f_\alpha(e^{-i\pi}t|\mu)\, dt 
\label{eq:pollardint2}
\intertext{
Hence $f_\alpha(x|\mu)$ is the Laplace transform of}
 \frac{1}{\pi}\,{\rm Im}\, \widetilde f_\alpha(e^{-i\pi}t|\mu) 
   &= \frac{1}{\pi}\,{\rm Im} \left\{ e^{-\mu(e^{-i\pi}t)^\alpha}\right\} \\
   &= \frac{1}{\pi}\, e^{-\mu t^\alpha \cos\pi\alpha} \sin(\mu t^\alpha \sin\pi\alpha) \\
 \frac{1}{\pi}\,{\rm Im}\, \widetilde f\,^\prime_\alpha(e^{-i\pi}t|\mu)
   &= -\frac{\mu\alpha}{\pi}\, {\rm Im}  \left\{ t^{\alpha-1}\, e^{-i\pi\alpha-\mu(e^{-i\pi}t)^\alpha}\right\} \\
 \lim_{n\to\infty} n\, \frac{1}{\pi}\,{\rm Im} \widetilde f\,^\prime_\alpha(e^{-i\pi}t|\tfrac{\mu}{n})
   &= -\frac{\mu\alpha}{\pi}\, {\rm Im}  \left\{ t^{\alpha-1}\, e^{-i\pi\alpha}\right\} 
  = \frac{\mu\alpha}{\pi}\, t^{\alpha-1}\sin\pi\alpha
\end{align}
($s\to e^{-i\pi}t=-t$ so that $ds\to -dt$).

We may thus extend the stable GGC LKCD of~(\ref{eq:3levelstableLKCD}) to include an additional column on the left,
as shown in~(\ref{eq:extended3level3columnstableLKCD}).
\begin{equation}
\begin{tikzpicture}[auto,scale=2.0, baseline=(current  bounding  box.center)]
\newcommand*{\size}{\scriptsize}%
\newcommand*{\gap}{.2ex}%
\newcommand*{\width}{2.25}%
\newcommand*{\height}{1.0}%

\node (O) at (0,0) {$\dfrac{1}{\pi}\,{\rm Im}\{\exp({-\mu e^{-i\pi\alpha}\,t^\alpha)}\}$};
\node (P) at ($(O)+(\width,0)$)  {$f_{\alpha}(x|\mu)$};
\node (Q) at ($(P)+(\width,0)$)  {$\exp(-\mu s^\alpha)$};
\node (A) at ($(O)-(0,\height)$) {$\dfrac{\mu\alpha}{\pi}\, \sin(\pi \alpha)\, t^{\alpha-1}$};
\node (B) at ($(A)+(\width,0)$)  {$\mu\alpha\, \dfrac{\,x^{-\alpha}}{\Gamma(1-\alpha)}$};
\node (C) at ($(B)+(\width,0)$) {$\mu\alpha s^{\alpha-1}$};
\node (X) at ($(A)-(0,\height)$)  {$(1-\alpha)\, \delta(t)$};
\node (Y) at ($(X)+(\width,0)$) {$1-\alpha$};
\node (Z) at  ($(Y)+(\width,0)$) {$\dfrac{1-\alpha}{s}$};   
\node (P0) at ($(P)$)  {$\phantom{\mu\alpha\, \dfrac{\,x^{-\alpha}}{\Gamma(1-\alpha)}}$};
\node (Y0) at ($(Y)$)  {$\phantom{\mu\alpha\, \dfrac{\,x^{-\alpha}}{\Gamma(1-\alpha)}}$};

\node (B1) at ($(B)+(0,7.5*\gap)$) {}; 
\node (C1) at ($(B1)+(\width,0)$) {}; 

\draw[Myarrow] ([yshift =  \gap]O.east)   --   ([yshift =  \gap]P0.west) ;
\draw[Myarrow] ([yshift = -\gap]P0.west)  --   ([yshift = -\gap]O.east); 
\draw[Myarrow] ([yshift =  \gap]P0.east)   --   ([yshift =  \gap]Q.west) ;
\draw[Myarrow] ([yshift = -\gap]Q.west)  --   ([yshift = -\gap]P0.east); 
\draw[Myarrow]([xshift  =  \gap]Q.south) --  ([xshift =  \gap]C1.north);
\draw[Myarrow]([xshift  = -\gap]C1.north) --   ([xshift = -\gap]Q.south);
\draw[Myarrow,dashed] ([xshift = -\gap]B1.north) -- ([xshift =  -\gap]P.south);
\draw[Myarrow] ([xshift =  \gap]P.south) -- ([xshift =  \gap]B1.north);
\draw[Myarrow] ([yshift =  \gap]A.east)   --   ([yshift =  \gap]B.west) ;
\draw[Myarrow] ([yshift = -\gap]B.west)  --   ([yshift = -\gap]A.east); 
\draw[Myarrow] ([yshift =  \gap]B.east)  --  ([yshift =  \gap]C.west) ;
\draw[Myarrow] ([yshift = -\gap]C.west) --   ([yshift = -\gap]B.east); 

\draw[Myarrow] ([xshift =  \gap]O.south) -- ([xshift =  \gap]A.north);

\draw[Myarrow]([xshift  =  \gap]C.south) --  ([xshift =  \gap]Z.north);
\draw[Myarrow]([xshift  = -\gap]Z.north) --   ([xshift = -\gap]C.south);
\draw[Myarrow,dashed] ([xshift = -\gap]Y.north) -- ([xshift =  -\gap]B.south);
\draw[Myarrow] ([xshift =  \gap]B.south) -- ([xshift =  \gap]Y.north);
\draw[Myarrow] ([yshift =  \gap]X.east)   --   ([yshift =  \gap]Y0.west) ;
\draw[Myarrow] ([yshift = -\gap]Y0.west)  --   ([yshift = -\gap]X.east); 
\draw[Myarrow] ([yshift =  \gap]Y0.east)  --  ([yshift =  \gap]Z.west) ;
\draw[Myarrow] ([yshift = -\gap]Z.west) --   ([yshift = -\gap]Y0.east); 

\fill[opacity=0.15]  ($(O)+(-\width*0.5,\height*0.2)$) rectangle ($(Y)+(-\width*0.24,-\height*0.22)$);
\fill[opacity=0.15]  ($(P)+(-\width*0.23,\height*0.2)$) rectangle ($(C)+(\width*0.24,\height*0.22)$);

\end{tikzpicture}
\label{eq:extended3level3columnstableLKCD}
\end{equation}
All horizontal arrows represent the Laplace transform or its inverse.
All down arrows represent a limiting process from the top layer of densities to the bottom layer,
or equivalently involve the logarithmic  derivative for the outer down arrows.

We note that  $x^{-\alpha}$ is the Laplace transform $\widetilde g_{\alpha,0}(x)$ of the density 
$g_{\alpha,0}(t)=  t^{\alpha-1}/\Gamma(\alpha)$.
That together with the  Euler identity $\Gamma(1-\alpha)\Gamma(\alpha) =  \pi/\sin(\pi\alpha)$ ensures consistency
between the  first two nodes in the middle layer.

The Dirac delta density in the bottom left node arises from the Laplace relation
\begin{align}
\int_0^\infty e^{-xt} \, \delta(t-a) dt &= e^{-xa}
\end{align}
We do not attempt to define a downarrow to the delta function from the node above.

For completeness, we note that the objects in the  first and last columns are related by two Laplace transforms via the middle column.
The relationship can be formulated as  composition of the two Laplace transforms.
Let the leftmost object be $u(t)$, so that  $f(x)=\widetilde u(x)$ and, in turn,  $\widetilde f(s)$ is 
\begin{align}
\widetilde f(s)  
                     &= \int_0^\infty dt \;e^{-sx}  \int_0^\infty e^{-xt} u(t) dx 
                       = \int_0^\infty u(t) \int_0^\infty e^{-(s+t)x} dx \; dt  \nonumber \\
                     &= \int_0^\infty \frac{u(t)}{s+t} \;dt 
\label{eq:stietjes}
\end{align}
Given $u(t)$, the Stieltjes transform,  
as (\ref{eq:stietjes})  is known, bypasses $f(x)$ to give $\widetilde f(s)$ directly.
In particular, the Laplace  transform of the representation~(\ref{eq:pollardint2}) of $f_\alpha(x|\mu)$ is
\begin{align}
\widetilde f_\alpha(s|\mu) 
  &= \frac{1}{\pi}\,{\rm Im} \int_0^\infty \frac{\widetilde f_\alpha(e^{-i\pi}t|\mu)}{s+t} \;dt 
\end{align}
But our prime objective  is to study the  ID density $f_\alpha(x|\mu)$ in the middle column.
It thus seems counter to that objective to seek to bypass $f_\alpha(x|\mu)$  as the Stieltjes transform does.

 We turn next to the convolution of two gamma densities to generate new ID densities. 
The  L\'{e}vy density of the convolution is the sum of the L\'{e}vy densities  of the individual ID densities.
 Sums are simple, but explicit convolutions can be complicated as the next example illustrates.

 \section{Convolution of Two Gamma Densities}
 \label{sec:convolution}

The Laplace convolution $\{f\star g\}(x)\equiv\{g\star f\}(x)$ of two functions $f(x)$ and $g(x)$ is defined by 
\begin{align}
\{f\star g\}(x) &= \int_0^x f(x-y)g(y)dy  = x\int_0^1 f(x(1-t))g(xt)dt 
\label{eq:convolution}
\end{align}
If $f$ and $g$ are ID densities,  so is $f\star g$.
The L\'{e}vy density of  $f\star g$ is the sum of the L\'{e}vy densities of $f$ and $g$.
The convolution of two gamma densities $g_{\mu_1,\lambda_1}$ and  $g_{\mu_2,\lambda_2}$
has the following LKCD
\begin{equation}
\begin{tikzpicture}[auto,scale=2.0, baseline=(current  bounding  box.center)]
\newcommand*{\size}{\scriptsize}%
\newcommand*{\gap}{.2ex}%
\newcommand*{\width}{2.75}%
\newcommand*{\height}{1.25}%

\node (P) at (0,0)  {$\{g_{\mu_1,\lambda_1}\star g_{\mu_2,\lambda_2}\}(x)$};
\node (Q) at ($(P)+(\width,0)$) {$\left(\dfrac{1}{\lambda_1+s}\right)^{\mu_1} \left(\dfrac{1}{\lambda_2+s}\right)^{\mu_2}$};
\node (B) at ($(P)-(0,\height)$) {$\mu_1 \, e^{-\lambda_1 x}+\mu_2 \, e^{-\lambda_2 x}$};
\node (C) at  ($(B)+(\width,0)$) {$\dfrac{\mu_1}{\lambda_1+s}+\dfrac{\mu_2}{\lambda_2+s}$};   

\draw[Myarrow] ([yshift =  \gap]P.east)   --   ([yshift =  \gap]Q.west) ;
\draw[Myarrow] ([yshift = -\gap]Q.west)  --   ([yshift = -\gap]P.east); 
\draw[Myarrow]([xshift  =  \gap]Q.south) --  ([xshift =  \gap]C.north);
\draw[Myarrow]([xshift  = -\gap]C.north) --   ([xshift = -\gap]Q.south);
\draw[Myarrow,dashed] ([xshift = -\gap]B.north) -- ([xshift =  -\gap]P.south);
\draw[Myarrow] ([xshift =  \gap]P.south) -- ([xshift =  \gap]B.north);
\draw[Myarrow] ([yshift =  \gap]B.east)  --  ([yshift =  \gap]C.west) ;
\draw[Myarrow] ([yshift = -\gap]C.west) --   ([yshift = -\gap]B.east); 
\end{tikzpicture}
\label{eq:gammaconvLKCD}
\end{equation}

 For $\lambda_1=\lambda_2=\lambda$, the gamma density
 is closed under convolution, {\it i.e.}\ $g_{\mu_1,\lambda}\star g_{\mu_2,\lambda}=g_{\mu_1+\mu_2,\lambda}$,
 as evident from~(\ref{eq:gammaconvLKCD}).
 However, for $\lambda_1\ne\lambda_2$, 
 convolution closure no longer holds, the resultant density is no longer gamma, although it remains ID. 
 The explicit form of $\{g_{\mu_1,\lambda_1}\star g_{\mu_2,\lambda_2}\}(x)$ is
 \begin{align}
 \{g_{\mu_1,\lambda_1}\star g_{\mu_2,\lambda_2}\}(x)
 &= x \int_0^1 g_{\mu_1,\lambda_1}(x(1-t)) g_{\mu_2,\lambda_2}(xt) dt \\
 &= \frac{x^{\mu_1+\mu_2-1}}{\Gamma(\mu_1)\Gamma(\mu_2)}\, e^{-\lambda_1 x}
  \int_0^1 (1-t)^{\mu_1-1} t^{\mu_2-1} \,  e^{-(\lambda_2-\lambda_1) xt} \, dt  \\
 &\equiv \frac{x^{\mu_1+\mu_2-1}}{\Gamma(\mu_1)\Gamma(\mu_2)}\, e^{-\lambda_2 x}
  \int_0^1 (1-t)^{\mu_2-1} t^{\mu_1-1} \,  e^{-(\lambda_1-\lambda_2) xt} \, dt  
\label{eq:gammaconv}
  \end{align}
We may  interpret the integral form~(\ref{eq:gammaconv}) in more than one way, as discussed next.

\subsection{Hypergeometric Function Interpretation}
\label{hypergeom}
Here and elsewhere,  we draw extensively on Abramowitz and Stegun~\cite{AbrSteg} for Laplace transform pairs 
and integrals such as the one in~(\ref{eq:gammaconv}), which may be  expressed as 
 \begin{align}
 \{g_{\mu_1,\lambda_1}\star g_{\mu_2,\lambda_2}\}(x)
 &= \frac{x^{\mu_1+\mu_2-1}}{\Gamma(\mu_1+\mu_2)}\,  e^{-\lambda_2 x} \,
 M(\mu_2,\mu_1+\mu_2,(\lambda_2-\lambda_1)\,x) \\
 &= 
 g_{\mu_1+\mu_2,\lambda_2}(x) \, M(\mu_2,\mu_1+\mu_2,(\lambda_2-\lambda_1)\,x) 
 \label{eq:hypergeom} \\
 {\rm where}\quad M(a,a+b,x) &= \frac{\Gamma(a+b)}{\Gamma(a)\Gamma(b)} \int_0^1 (1-t)^{a-1} t^{b-1} \,  e^{xt} \, dt   
\end{align}
is the confluent hypergeometric function  (\cite{AbrSteg}~p505,~13.2.1).
Since $M(\cdot,\cdot,0)=1$,  $\lambda_1=\lambda_2=\lambda$ reproduces gamma closure 
$ \{g_{\mu_1,\lambda}\star g_{\mu_2,\lambda}\}(x)=g_{\mu_1+\mu_2,\lambda}(x)$. 
Hypergeometric functions are flexible infinite series representations of  a variety of functions 
for different  choices of arguments.
In particular, with the aid of~\cite{AbrSteg}~p509, 13.6.3 and the Legendre duplication formula 
$\sqrt{\pi}\,\Gamma(2\mu)=2^{2\mu-1}\Gamma(\mu)\Gamma(\mu+1/2)$, 
the case $\mu_1=\mu_2=\mu$ ($\lambda_1\ne\lambda_2$) can  be shown to be
\begin{align}
 \{g_{\mu,\lambda_1}\star g_{\mu,\lambda_2}\}(x) &= 
\frac{\sqrt{\pi}}{\Gamma(\mu)} \,  e^{-(\lambda_1+\lambda_2)x/2}\, 
 \left(\frac{x}{\lambda_2-\lambda_1}\right)^{\mu-\tfrac{1}{2}}\, I_{\mu-\tfrac{1}{2}}\left(\frac{\lambda_2-\lambda_1}{2}\, x\right)
\label{eq:gammaconvBessel}
 \end{align}
$I_\mu(x)$ is the modified Bessel function of the first kind of order $\mu$.
Alternatively, (\ref{eq:gammaconvBessel}) is given by  the Laplace transform pair~\cite{AbrSteg}~p1024~(29.3.50).
It leads directly  to the sum of two 
exponentials  under the limit~(\ref{eq:limnxf}), noting  that $\sqrt{\pi x/2}\,I_{-1/2}(x)=\cosh(x)$:
\begin{align}
\lim_{n\to\infty}n\, x\, \{g_{\mu/n,\lambda_1}\star g_{\mu/n,\lambda_2}\}(x)  
&= \mu\,   e^{-(\lambda_1+\lambda_2)x/2}\, \sqrt{\pi(\lambda_2-\lambda_1)x}\;
I_{-\tfrac{1}{2}}\left(\frac{\lambda_2-\lambda_1}{2}\, x\right) \\
 &= \mu\,e^{-\lambda_1 x}+\mu\, e^{-\lambda_2 x}
\label{eq:limbessel}
\end{align} 
Let $(\lambda_1,\lambda_2)=(0,1)$,  for which  the convolution~(\ref{eq:gammaconvBessel}) is
\vspace{-1pc}
 \begin{labeling}{$\mu=\tfrac{1}{2}$:}
\setlength{\partopsep}{0pt}
\setlength{\topsep}{0pt}
\setlength{\itemsep}{0pt}
\setlength{\parskip}{0pt}
\setlength{\parsep}{0pt}
\item[$\mu=1$:] 
  $e^{-x/2}\,\sqrt{\pi x}\,I_{1/2}(x/2)\equiv \{g_{1,0}\star g_{1,1}\}(x)=1-e^{-x}$, 
   where the latter expression follows from $\sqrt{\pi x/2}\,I_{1/2}(x)=\sinh(x)$.
   It is the convolution of a constant and an exponential.
   It also  follows from following the Laplace transform route from the bottom left node of~(\ref{eq:gammaconvLKCD}).
\item[$\mu=\tfrac{1}{2}$:] so that (\ref{eq:gammaconvBessel}) becomes $e^{-x/2}I_0(x/2) \equiv \{g_{1/2,0}\star g_{1/2,1}\}(x)$. 
\end{labeling}
\vspace{-1pc}
The LKCD of these two  cases is shown below, with the  leftmost limiting downarrow omitted
since $\mu$ is fixed.
$1\pm e^{-x}$ may also be written as $e^{-x/2}\sqrt{\pi x}\, I_{\pm 1/2}(x/2)$.
\begin{equation}
\begin{tikzpicture}[auto,scale=2.0, baseline=(current  bounding  box.center)]
\newcommand*{\size}{\scriptsize}%
\newcommand*{\gap}{.2ex}%
\newcommand*{\width}{2.0}%
\newcommand*{\height}{1.25}%

\node (P) at (0,0)  {$1- e^{-x}$};
\node (Q) at ($(P)+(\width,0)$) {$\dfrac{1}{s}-\dfrac{1}{1+s}$};
\node (B) at ($(P)-(0,\height)$) {$1+ e^{-x}$};
\node (C) at  ($(B)+(\width,0)$) {$\dfrac{1}{s}+\dfrac{1}{1+s}$};   

\draw[Myarrow] ([yshift =  \gap]P.east)   --   ([yshift =  \gap]Q.west) ;
\draw[Myarrow] ([yshift = -\gap]Q.west)  --   ([yshift = -\gap]P.east); 
\draw[Myarrow]([xshift  =  \gap]Q.south) --  ([xshift =  \gap]C.north);
\draw[Myarrow]([xshift  = -\gap]C.north) --   ([xshift = -\gap]Q.south);
\draw[Myarrow] ([yshift =  \gap]B.east)  --  ([yshift =  \gap]C.west) ;
\draw[Myarrow] ([yshift = -\gap]C.west) --   ([yshift = -\gap]B.east); 

\node (P) at (\width*1.8,0)  {$e^{-x/2}I_0(\tfrac{x}{2})$};
\node (Q) at ($(P)+(\width,0)$) {$\dfrac{1}{\sqrt{s(1+s)}}$};
\node (B) at ($(P)-(0,\height)$) {$\dfrac{1}{2} \, (1+ e^{-x})$};
\node (C) at  ($(B)+(\width,0)$) {$\dfrac{1}{2}\left(\dfrac{1}{s}+\dfrac{1}{1+s}\right)$};   

\draw[Myarrow] ([yshift =  \gap]P.east)   --   ([yshift =  \gap]Q.west) ;
\draw[Myarrow] ([yshift = -\gap]Q.west)  --   ([yshift = -\gap]P.east); 
\draw[Myarrow]([xshift  =  \gap]Q.south) --  ([xshift =  \gap]C.north);
\draw[Myarrow]([xshift  = -\gap]C.north) --   ([xshift = -\gap]Q.south);
\draw[Myarrow] ([yshift =  \gap]B.east)  --  ([yshift =  \gap]C.west) ;
\draw[Myarrow] ([yshift = -\gap]C.west) --   ([yshift = -\gap]B.east); 
\end{tikzpicture}
\label{eq:gammahalfconvLKCD}
\end{equation}

The mathematical generality of the confluent  hypergeometric function can hide probabilistic insight.
We can, in fact, interpret it  as the Laplace transform of a familiar density, as we do next  in the second of our two interpretations 
of the convolution of two gamma densities.

\subsection{Beta Density Interpretation}
\label{beta}
The beta distribution with two shape parameters $a,b$ has the density
 \begin{align}
\mathrm{Beta}(x|a,b) &= \frac{\Gamma(a+b)}{\Gamma(a)\Gamma(b)} x^{a-1}(1-x)^{b-1}  \qquad 0\le x \le 1
\label{eq:beta} \\
\mathrm{with}\quad 
\mathrm{Beta}(x|1-a,a) &= \frac{\sin(\pi a)}{\pi} x^{-a}(1-x)^{a-1} 
\intertext{We may extend the  domain to $\mathbb{R}_+$ by defining}
 \beta_{a,b}(x) 
 &= \begin{cases}
         \mathrm{Beta}(x|a,b)  &  0\le x \le 1 \\
         0 & x>1
       \end{cases}                      
\label{eq:extendedbeta}
 \end{align}

The Laplace transform of~(\ref{eq:extendedbeta})  is
 \begin{align}
 \widetilde\beta_{a,b}(x) &= \int_0^\infty e^{-xt} \beta_{a,b}(t) dt = M(b,a+b,-x) 
 \label{eq:extendedbetaLT} \\
 \implies  
 \{g_{\mu_1,\lambda_1}\star g_{\mu_2,\lambda_2}\}(x)
 &= 
 g_{\mu_1+\mu_2,\lambda_2}(x) \, \widetilde\beta_{\mu_1,\mu_2}((\lambda_1-\lambda_2)\,x)
\label{eq:gammaconvBeta1} \\
 &\equiv g_{\mu_1+\mu_2,\lambda_1}(x) \, \widetilde\beta_{\mu_2,\mu_1}((\lambda_2-\lambda_1)\,x)
\label{eq:gammaconvBeta2}
\intertext{For $\mu_1+\mu_2=1$, let  $\mu_1\equiv 1-\alpha, \mu_2\equiv\alpha$ ($0<\alpha<1$). 
Further, set $\lambda_1=0,\lambda_2=1$, to get}
  \{g_{1-\alpha,0}\star g_{\alpha,1}\}(x) &= \widetilde\beta_{\alpha,1-\alpha}(x).
\label{eq:gammaconvBeta3}
 \end{align}
We recall that, for the stable case represented in~(\ref{eq:extended3level3columnstableLKCD}), $\rho_\alpha(x)=\alpha g_{1-\alpha,0}(x)$.
Hence we can  take the convolution with  $g_{\alpha,1}(x)$ to be at the centre of~(\ref{eq:extended3level3columnstableLKCD}), 
so that we expand from that point  in all four directions.
The LKCD thus has the following form (explicit forms  for the top-level densities are not readily available).
\begin{equation}
\begin{tikzpicture}[auto,scale=2.0, baseline=(current  bounding  box.center)]
\newcommand*{\size}{\scriptsize}%
\newcommand*{\gap}{.2ex}%
\newcommand*{\width}{2.25}%
\newcommand*{\height}{1.0}%

\node (O) at (0,0) {$\dfrac{1}{\pi}\mathrm{Im}\left\{\widetilde h_\alpha(e^{-i\pi}t|\mu)\right\}$};
\node (P) at ($(O)+(\width,0)$)  {$h_\alpha(x|\mu)$};
\node (Q) at ($(P)+(\width,0)$)  {$\widetilde h_\alpha(s|\mu)$};
\node (A) at ($(O)-(0,\height)$) {$\mu\beta_{\alpha,1-\alpha}(t)$};
\node (B) at ($(A)+(\width,0)$)  {$\mu\{g_{1-\alpha,0}\star g_{\alpha,1}\}(x)$};
\node (C) at ($(B)+(\width,0)$) {$\mu\dfrac{s^{\alpha-1}}{(1+s)^{\alpha}}$};
\node (X) at ($(A)-(0,\height)$)  {$(1-\alpha)\delta(t) +\alpha \delta(t-1)$};
\node (Y) at ($(X)+(\width,0)$) {$1-\alpha +\alpha \, e^{-x}$};
\node (Z) at  ($(Y)+(\width,0)$) {$\dfrac{1-\alpha}{s}+\dfrac{\alpha}{1+s}$};   

\draw[Myarrow] ([yshift =  \gap]O.east)   --   ([yshift =  \gap]P.west) ;
\draw[Myarrow] ([yshift = -\gap]P.west)  --   ([yshift = -\gap]O.east); 
\draw[Myarrow] ([yshift =  \gap]P.east)   --   ([yshift =  \gap]Q.west) ;
\draw[Myarrow] ([yshift = -\gap]Q.west)  --   ([yshift = -\gap]P.east); 
\draw[Myarrow]([xshift  =  \gap]Q.south) --  ([xshift =  \gap]C.north);
\draw[Myarrow]([xshift  = -\gap]C.north) --   ([xshift = -\gap]Q.south);
\draw[Myarrow,dashed] ([xshift = -\gap]B.north) -- ([xshift =  -\gap]P.south);
\draw[Myarrow] ([xshift =  \gap]P.south) -- ([xshift =  \gap]B.north);
\draw[Myarrow] ([yshift =  \gap]A.east)   --   ([yshift =  \gap]B.west) ;
\draw[Myarrow] ([yshift = -\gap]B.west)  --   ([yshift = -\gap]A.east); 
\draw[Myarrow] ([yshift =  \gap]B.east)  --  ([yshift =  \gap]C.west) ;
\draw[Myarrow] ([yshift = -\gap]C.west) --   ([yshift = -\gap]B.east); 

\draw[Myarrow] ([xshift =  \gap]O.south) -- ([xshift =  \gap]A.north);

\draw[Myarrow]([xshift  =  \gap]C.south) --  ([xshift =  \gap]Z.north);
\draw[Myarrow]([xshift  = -\gap]Z.north) --   ([xshift = -\gap]C.south);
\draw[Myarrow,dashed] ([xshift = -\gap]Y.north) -- ([xshift =  -\gap]B.south);
\draw[Myarrow] ([xshift =  \gap]B.south) -- ([xshift =  \gap]Y.north);
\draw[Myarrow] ([yshift =  \gap]X.east)   --   ([yshift =  \gap]Y.west) ;
\draw[Myarrow] ([yshift = -\gap]Y.west)  --   ([yshift = -\gap]X.east); 
\draw[Myarrow] ([yshift =  \gap]Y.east)  --  ([yshift =  \gap]Z.west) ;
\draw[Myarrow] ([yshift = -\gap]Z.west) --   ([yshift = -\gap]Y.east); 
\end{tikzpicture}
\label{eq:gammaconv1LKCD}
\end{equation}
 In keeping with the foregoing discussion, we may also write $\beta_{\alpha,1-\alpha}(t)$ as
 \begin{align}
 \beta_{\alpha,1-\alpha}(t)
   &= \frac{1}{\pi}\mathrm{Im}\left\{\frac{(e^{-i\pi}t)^{\alpha-1}}{(1+e^{-i\pi}t)^{\alpha}}\right\}
   = \begin{cases}
        \frac{\sin\pi\alpha}{\pi}\, t^{\alpha-1}(1-t)^{-\alpha}  &  0\le t \le 1 \\
         0 & t>1
       \end{cases}                 
 \end{align}
 
 Setting $(\lambda_1,\lambda_2)=(1,0)$  gives (where $h_\alpha(x|\mu)$ denotes a correspondingly different density)
\begin{equation}
\begin{tikzpicture}[auto,scale=2.0, baseline=(current  bounding  box.center)]
\newcommand*{\size}{\scriptsize}%
\newcommand*{\gap}{.2ex}%
\newcommand*{\width}{2.25}%
\newcommand*{\height}{1.0}%

\node (O) at (0,0) {$\dfrac{1}{\pi}\mathrm{Im}\left\{\widetilde h_\alpha(e^{-i\pi}t|\mu)\right\}$};
\node (P) at ($(O)+(\width,0)$)  {$h_\alpha(x|\mu)$};
\node (Q) at ($(P)+(\width,0)$)  {$\widetilde h_\alpha(s|\mu)$};
\node (A) at ($(O)-(0,\height)$) {$\mu\beta_{1-\alpha,\alpha}(t)$};
\node (B) at ($(A)+(\width,0)$)  {$\mu\{g_{1-\alpha,1}\star g_{\alpha,0}\}(x)$};
\node (C) at ($(B)+(\width,0)$) {$\mu\dfrac{(1+s)^{\alpha-1}}{s^{\alpha}}$};
\node (X) at ($(A)-(0,\height)$) {$\alpha \delta(t) + (1-\alpha) \delta(t-1)$};
\node (Y) at ($(X)+(\width,0)$) {$\alpha +(1-\alpha) \, e^{-x}$};
\node (Z) at  ($(Y)+(\width,0)$) {$\dfrac{\alpha}{s}+\dfrac{1-\alpha}{1+s}$};   

\draw[Myarrow] ([yshift =  \gap]O.east)   --   ([yshift =  \gap]P.west) ;
\draw[Myarrow] ([yshift = -\gap]P.west)  --   ([yshift = -\gap]O.east); 
\draw[Myarrow] ([yshift =  \gap]P.east)   --   ([yshift =  \gap]Q.west) ;
\draw[Myarrow] ([yshift = -\gap]Q.west)  --   ([yshift = -\gap]P.east); 
\draw[Myarrow]([xshift  =  \gap]Q.south) --  ([xshift =  \gap]C.north);
\draw[Myarrow]([xshift  = -\gap]C.north) --   ([xshift = -\gap]Q.south);
\draw[Myarrow,dashed] ([xshift = -\gap]B.north) -- ([xshift =  -\gap]P.south);
\draw[Myarrow] ([xshift =  \gap]P.south) -- ([xshift =  \gap]B.north);
\draw[Myarrow] ([yshift =  \gap]A.east)   --   ([yshift =  \gap]B.west) ;
\draw[Myarrow] ([yshift = -\gap]B.west)  --   ([yshift = -\gap]A.east); 
\draw[Myarrow] ([yshift =  \gap]B.east)  --  ([yshift =  \gap]C.west) ;
\draw[Myarrow] ([yshift = -\gap]C.west) --   ([yshift = -\gap]B.east); 

\draw[Myarrow] ([xshift =  \gap]O.south) -- ([xshift =  \gap]A.north);

\draw[Myarrow]([xshift  =  \gap]C.south) --  ([xshift =  \gap]Z.north);
\draw[Myarrow]([xshift  = -\gap]Z.north) --   ([xshift = -\gap]C.south);
\draw[Myarrow,dashed] ([xshift = -\gap]Y.north) -- ([xshift =  -\gap]B.south);
\draw[Myarrow] ([xshift =  \gap]B.south) -- ([xshift =  \gap]Y.north);
\draw[Myarrow] ([yshift =  \gap]X.east)   --   ([yshift =  \gap]Y.west) ;
\draw[Myarrow] ([yshift = -\gap]Y.west)  --   ([yshift = -\gap]X.east); 
\draw[Myarrow] ([yshift =  \gap]Y.east)  --  ([yshift =  \gap]Z.west) ;
\draw[Myarrow] ([yshift = -\gap]Z.west) --   ([yshift = -\gap]Y.east); 
\end{tikzpicture}
\label{eq:gammaconv2LKCD}
\end{equation}

For $\alpha=1/2$, we have explicit forms, recalling that 
$e^{-x/2}I_0(x/2) \equiv \{g_{1/2,0}\star g_{1/2,1}\}(x)$:
\begin{equation}
\begin{tikzpicture}[auto,scale=2.0, baseline=(current  bounding  box.center)]
\newcommand*{\size}{\scriptsize}%
\newcommand*{\gap}{.2ex}%
\newcommand*{\width}{2.25}%
\newcommand*{\height}{1.0}%

\node (O) at (0,0) {$\dfrac{1}{\pi}\mathrm{Im} \left(\sqrt{1-t}+i\sqrt{t}\right)^{2\mu}$};
\node (P) at ($(O)+(\width,0)$) {$\dfrac{\mu}{x}e^{-x/2}\,I_\mu\left(\tfrac{x}{2}\right)$};
\node (Q) at ($(P)+(\width,0)$) {$\left(\sqrt{1+s}-\sqrt{s}\right)^{2\mu}$};
\node (A) at ($(O)-(0,\height)$) {$\mu\beta_{1/2,1/2}(t)$};
\node (B) at ($(A)+(\width,0)$)  {$\mu e^{-x/2}I_0(\tfrac{x}{2})$};
\node (C) at ($(B)+(\width,0)$) {$\mu \dfrac{1}{\sqrt{s(1+s)}}$};
\node (X) at ($(A)-(0,\height)$)  {$\dfrac{1}{2}\left(\delta(t) +\delta(t-1)\right)$};
\node (Y) at ($(X)+(\width,0)$) {$\dfrac{1}{2}\left(1+ e^{-x}\right)$};
\node (Z) at  ($(Y)+(\width,0)$) {$\dfrac{1}{2}\left(\dfrac{1}{s}+\dfrac{1}{1+s}\right)$};   

\draw[Myarrow] ([yshift =  \gap]O.east)   --   ([yshift =  \gap]P.west) ;
\draw[Myarrow] ([yshift = -\gap]P.west)  --   ([yshift = -\gap]O.east); 
\draw[Myarrow] ([yshift =  \gap]P.east)   --   ([yshift =  \gap]Q.west) ;
\draw[Myarrow] ([yshift = -\gap]Q.west)  --   ([yshift = -\gap]P.east); 
\draw[Myarrow]([xshift  =  \gap]Q.south) --  ([xshift =  \gap]C.north);
\draw[Myarrow]([xshift  = -\gap]C.north) --   ([xshift = -\gap]Q.south);
\draw[Myarrow,dashed] ([xshift = -\gap]B.north) -- ([xshift =  -\gap]P.south);
\draw[Myarrow] ([xshift =  \gap]P.south) -- ([xshift =  \gap]B.north);
\draw[Myarrow] ([yshift =  \gap]A.east)   --   ([yshift =  \gap]B.west) ;
\draw[Myarrow] ([yshift = -\gap]B.west)  --   ([yshift = -\gap]A.east); 
\draw[Myarrow] ([yshift =  \gap]B.east)  --  ([yshift =  \gap]C.west) ;
\draw[Myarrow] ([yshift = -\gap]C.west) --   ([yshift = -\gap]B.east); 

\draw[Myarrow] ([xshift =  \gap]O.south) -- ([xshift =  \gap]A.north);

\draw[Myarrow]([xshift  =  \gap]C.south) --  ([xshift =  \gap]Z.north);
\draw[Myarrow]([xshift  = -\gap]Z.north) --   ([xshift = -\gap]C.south);
\draw[Myarrow] ([yshift =  \gap]X.east)   --   ([yshift =  \gap]Y.west) ;
\draw[Myarrow] ([yshift = -\gap]Y.west)  --   ([yshift = -\gap]X.east); 
\draw[Myarrow] ([yshift =  \gap]Y.east)  --  ([yshift =  \gap]Z.west) ;
\draw[Myarrow] ([yshift = -\gap]Z.west) --   ([yshift = -\gap]Y.east); 
\end{tikzpicture}
\label{eq:gammaconvLKCD_half}
\end{equation}
We have used the  Laplace transform  pair~\cite{AbrSteg}~p1024~(29.3.53) for the density $\mu e^{-x/2}I_\mu(x/2)/x$  and
its Laplace transform in the top right node, which may  be written  in several  forms
 \begin{align*}
\left(\sqrt{1+s}+\sqrt{s}\right)^{-2\mu} 
= \left(\sqrt{1+s}-\sqrt{s}\right)^{2\mu} 
&= \left(1+2s+2\sqrt{s(1+s)}\right)^{-\mu}  
= \left(1+2s-2\sqrt{s(1+s)}\right)^{\mu} 
 \end{align*}

Bessel functions play a prominent role in random walks and Brownian motion, with numerous applications to
random phenomena in physics, chemistry, biology, finance~{\it etc.}
For example,  in a section titled  \textquote{Bessel Functions and Random Walks} (\cite{Feller2}~p58-61), 
Feller showed that  the distribution of the first passage through $\mu>0$ 
 ({\it i.e.}\ the time  it takes to reach the point~$\mu$ for the first time in a random walk in one dimension, 
 starting at the point~0)
has the density  $\mu\,e^{-x}I_\mu(x)/x$ (where $x$ is time in this context).  
Feller proceeded to calculate  the Laplace transform of this density (\cite{Feller2}~p437) and to demonstrate its infinite divisibility, 
with $\rho(x)= e^{-x}I_0(x)$  (\cite{Feller2}~p451).
Without explicit reference to the Bessel function form, Bondesson (\cite{Bondesson}~p37)  showed that the first passage distribution
is GGC with Thorin density $\beta_{1/2,1/2}(t)$.
The latter is often referred to as the density of the arc-sine distribution  
in the  literature on random walks and Brownian motion.

Lastly, let $\mu_1+\mu_2=1-\theta$ for some parameter $\theta$.
Once again set $\mu_1=1-\alpha$ so that $\mu_2=\alpha-\theta$.
Since the gamma shape parameters $\mu_1,\mu_2$ are always positive, 
we must have $0\le \theta<\alpha$ where $0<\alpha<1$ as before.
Then~(\ref{eq:gammaconvBeta1}) becomes
 \begin{align}
 \{g_{1-\alpha,0}\star g_{\alpha-\theta,1}\}(x)
 &=  g_{1-\theta,0}(x) \, \widetilde\beta_{\alpha-\theta,1-\alpha}(x) 
\label{eq:gammaconvalphatheta}
 \end{align}
 which reduces to~(\ref{eq:gammaconvBeta3}) for $\theta=0$.
For $0<\theta<\alpha$,~(\ref{eq:gammaconvalphatheta})  is, in turn, the Laplace transform of a convolution of a gamma and a beta density
\begin{equation}
  \frac{\sin\pi\theta}{\pi}\, t^{\theta-1}\star \beta_{\alpha-\theta,1-\alpha}(t)
  \longrightarrow g_{1-\theta,0}(x) \, \widetilde\beta_{\alpha-\theta,1-\alpha}(x)
\end{equation}
This generalises~(\ref{eq:gammaconv1LKCD}) ({\it i.e.}\ the case $\theta=0$) to the following: 
\begin{equation}
\begin{tikzpicture}[auto,scale=2.0, baseline=(current  bounding  box.center)]
\newcommand*{\size}{\scriptsize}%
\newcommand*{\gap}{.2ex}%
\newcommand*{\width}{2.5}%
\newcommand*{\height}{1.0}%

\node (O) at (0,0) {$\dfrac{1}{\pi}\mathrm{Im}\left\{\widetilde h_{\alpha,\theta}(e^{-i\pi}t|\mu)\right\}$};
\node (P) at ($(O)+(\width,0)$)  {$h_{\alpha,\theta}(x|\mu)$};
\node (Q) at ($(P)+(\width,0)$)  {$\widetilde h_{\alpha,\theta}(s|\mu)$};
\node (A) at ($(O)-(0,\height)$) {$\mu\frac{\sin\pi\theta}{\pi}\, t^{\theta-1}\star \beta_{\alpha-\theta,1-\alpha}(t)$};
\node (B) at ($(A)+(\width,0)$)  {$\mu\{g_{1-\alpha,0}\star g_{\alpha-\theta,1}\}(x)$};
\node (C) at ($(B)+(\width,0)$) {$\mu\dfrac{s^{\alpha-1}}{(1+s)^{\alpha-\theta}}$};
\node (X) at ($(A)-(0,\height)$)  {$(1-\alpha)\delta(t) +(\alpha-\theta) \delta(t-1)$};
\node (Y) at ($(X)+(\width,0)$) {$1-\alpha +(\alpha-\theta) \, e^{-x}$};
\node (Z) at  ($(Y)+(\width,0)$) {$\dfrac{1-\alpha}{s}+\dfrac{\alpha-\theta}{1+s}$};   

\draw[Myarrow] ([yshift =  \gap]O.east)   --   ([yshift =  \gap]P.west) ;
\draw[Myarrow] ([yshift = -\gap]P.west)  --   ([yshift = -\gap]O.east); 
\draw[Myarrow] ([yshift =  \gap]P.east)   --   ([yshift =  \gap]Q.west) ;
\draw[Myarrow] ([yshift = -\gap]Q.west)  --   ([yshift = -\gap]P.east); 
\draw[Myarrow]([xshift  =  \gap]Q.south) --  ([xshift =  \gap]C.north);
\draw[Myarrow]([xshift  = -\gap]C.north) --   ([xshift = -\gap]Q.south);
\draw[Myarrow,dashed] ([xshift = -\gap]B.north) -- ([xshift =  -\gap]P.south);
\draw[Myarrow] ([xshift =  \gap]P.south) -- ([xshift =  \gap]B.north);
\draw[Myarrow] ([yshift =  \gap]A.east)   --   ([yshift =  \gap]B.west) ;
\draw[Myarrow] ([yshift = -\gap]B.west)  --   ([yshift = -\gap]A.east); 
\draw[Myarrow] ([yshift =  \gap]B.east)  --  ([yshift =  \gap]C.west) ;
\draw[Myarrow] ([yshift = -\gap]C.west) --   ([yshift = -\gap]B.east); 

\draw[Myarrow] ([xshift =  \gap]O.south) -- ([xshift =  \gap]A.north);

\draw[Myarrow]([xshift  =  \gap]C.south) --  ([xshift =  \gap]Z.north);
\draw[Myarrow]([xshift  = -\gap]Z.north) --   ([xshift = -\gap]C.south);
\draw[Myarrow,dashed] ([xshift = -\gap]Y.north) -- ([xshift =  -\gap]B.south);
\draw[Myarrow] ([xshift =  \gap]B.south) -- ([xshift =  \gap]Y.north);
\draw[Myarrow] ([yshift =  \gap]X.east)   --   ([yshift =  \gap]Y.west) ;
\draw[Myarrow] ([yshift = -\gap]Y.west)  --   ([yshift = -\gap]X.east); 
\draw[Myarrow] ([yshift =  \gap]Y.east)  --  ([yshift =  \gap]Z.west) ;
\draw[Myarrow] ([yshift = -\gap]Z.west) --   ([yshift = -\gap]Y.east); 
\end{tikzpicture}
\label{eq:gammaconvalphathetaLKCD}
\end{equation}
In this case, $\alpha=1/2$ gives
\begin{equation}
\begin{tikzpicture}[auto,scale=2.0, baseline=(current  bounding  box.center)]
\newcommand*{\size}{\scriptsize}%
\newcommand*{\gap}{.2ex}%
\newcommand*{\width}{2.5}%
\newcommand*{\height}{1.0}%

\node (O) at (0,0) {$\dfrac{1}{\pi}\mathrm{Im}\left\{\widetilde h_{1/2,\theta}(e^{-i\pi}t|\mu)\right\}$};
\node (P) at ($(O)+(\width,0)$)  {$h_{1/2,\theta}(x|\mu)$};
\node (Q) at ($(P)+(\width,0)$)  {$\widetilde h_{1/2,\theta}(s|\mu)$};
\node (A) at ($(O)-(0,\height)$) {$\mu\frac{\sin\pi\theta}{\pi}\, t^{\theta-1}\star \beta_{1/2-\theta,1/2}(t)$};
\node (B) at ($(A)+(\width,0)$)  {$\mu\{g_{1/2,0}\star g_{1/2-\theta,1}\}(x)$};
\node (C) at ($(B)+(\width,0)$) {$\mu\dfrac{(1+s)^\theta}{\sqrt{s(1+s)}}$};
\node (X) at ($(A)-(0,\height)$)  {$\frac{1}{2}\delta(t) + \left(\frac{1}{2} - \theta\right)\delta(t-1)$};
\node (Y) at ($(X)+(\width,0)$) {$\frac{1}{2}+\left(\frac{1}{2} -\theta\right)e^{-x}$};
\node (Z) at  ($(Y)+(\width,0)$) {$\dfrac{1}{2s}+\dfrac{1/2 -\theta}{1+s}$};   

\draw[Myarrow] ([yshift =  \gap]O.east)   --   ([yshift =  \gap]P.west) ;
\draw[Myarrow] ([yshift = -\gap]P.west)  --   ([yshift = -\gap]O.east); 
\draw[Myarrow] ([yshift =  \gap]P.east)   --   ([yshift =  \gap]Q.west) ;
\draw[Myarrow] ([yshift = -\gap]Q.west)  --   ([yshift = -\gap]P.east); 
\draw[Myarrow]([xshift  =  \gap]Q.south) --  ([xshift =  \gap]C.north);
\draw[Myarrow]([xshift  = -\gap]C.north) --   ([xshift = -\gap]Q.south);
\draw[Myarrow,dashed] ([xshift = -\gap]B.north) -- ([xshift =  -\gap]P.south);
\draw[Myarrow] ([xshift =  \gap]P.south) -- ([xshift =  \gap]B.north);
\draw[Myarrow] ([yshift =  \gap]A.east)   --   ([yshift =  \gap]B.west) ;
\draw[Myarrow] ([yshift = -\gap]B.west)  --   ([yshift = -\gap]A.east); 
\draw[Myarrow] ([yshift =  \gap]B.east)  --  ([yshift =  \gap]C.west) ;
\draw[Myarrow] ([yshift = -\gap]C.west) --   ([yshift = -\gap]B.east); 

\draw[Myarrow] ([xshift =  \gap]O.south) -- ([xshift =  \gap]A.north);

\draw[Myarrow]([xshift  =  \gap]C.south) --  ([xshift =  \gap]Z.north);
\draw[Myarrow]([xshift  = -\gap]Z.north) --   ([xshift = -\gap]C.south);
\draw[Myarrow] ([yshift =  \gap]X.east)   --   ([yshift =  \gap]Y.west) ;
\draw[Myarrow] ([yshift = -\gap]Y.west)  --   ([yshift = -\gap]X.east); 
\draw[Myarrow] ([yshift =  \gap]Y.east)  --  ([yshift =  \gap]Z.west) ;
\draw[Myarrow] ([yshift = -\gap]Z.west) --   ([yshift = -\gap]Y.east); 
\end{tikzpicture}
\label{eq:gammaconvalfahalfthetaLKCD}
\end{equation}
$h_{\alpha,\theta}(x|\mu)$ is a three-parameter ID density generated by
$\{g_{1-\alpha,0}\star g_{\alpha-\theta,1}\}(x)$ which, in turn,  is a two-parameter ID density.
Aside from $h_{1/2,0}(x|\mu) = \mu e^{-x/2}I_\mu(x/2)/x$ of the LKCD given in~(\ref{eq:gammaconvLKCD_half}),
we shall not pursue further here the explicit  form for  the general case $h_{\alpha,\theta}(x|\mu)$.
We note only that, as shown in~(\ref{eq:gammaconvalphathetaLKCD}), it satisfies 
\begin{align}
-\frac{\widetilde h_{\alpha,\theta}^\prime(s|\mu)}{\widetilde h_{\alpha,\theta}(s|\mu)}
  &= \mu\dfrac{s^{\alpha-1}}{(1+s)^{\alpha-\theta}}
  \equiv \mu\,\rho_{\alpha,\theta}(s) \\
  \implies\quad
  \widetilde h_{\alpha,\theta}(s|\mu)
  &= \exp\left\{-\mu\int_0^s \rho_{\alpha,\theta}(t) dt \right\} 
\label{eq:h3param}
 \end{align}

Setting $(\lambda_1,\lambda_2)=(1,0)$ gives 
\begin{equation}
\begin{tikzpicture}[auto,scale=2.0, baseline=(current  bounding  box.center)]
\newcommand*{\size}{\scriptsize}%
\newcommand*{\gap}{.2ex}%
\newcommand*{\width}{2.5}%
\newcommand*{\height}{1.0}%

\node (O) at (0,0) {$\dfrac{1}{\pi}\mathrm{Im}\left\{\widetilde h_{\alpha,\theta}(e^{-i\pi}t|\mu)\right\}$};
\node (P) at ($(O)+(\width,0)$)  {$h_{\alpha,\theta}(x|\mu)$};
\node (Q) at ($(P)+(\width,0)$)  {$\widetilde h_{\alpha,\theta}(s|\mu)$};
\node (A) at ($(O)-(0,\height)$) {$\mu\frac{\sin\pi\theta}{\pi}\, t^{\theta-1}\star \beta_{1-\alpha,\alpha-\theta}(t)$};
\node (B) at ($(A)+(\width,0)$)  {$\mu\{g_{1-\alpha,1}\star g_{\alpha-\theta,0}\}(x)$};
\node (C) at ($(B)+(\width,0)$) {$\mu\dfrac{(1+s)^{\alpha-1}}{s^{\alpha-\theta}}$};
\node (X) at ($(A)-(0,\height)$)  {$(\alpha-\theta)\delta(t) +(1-\alpha) \delta(t-1)$};
\node (Y) at ($(X)+(\width,0)$) {$\alpha-\theta +(1-\alpha) \, e^{-x}$};
\node (Z) at  ($(Y)+(\width,0)$) {$\dfrac{\alpha-\theta}{s}+\dfrac{1-\alpha}{1+s}$};   

\draw[Myarrow] ([yshift =  \gap]O.east)   --   ([yshift =  \gap]P.west) ;
\draw[Myarrow] ([yshift = -\gap]P.west)  --   ([yshift = -\gap]O.east); 
\draw[Myarrow] ([yshift =  \gap]P.east)   --   ([yshift =  \gap]Q.west) ;
\draw[Myarrow] ([yshift = -\gap]Q.west)  --   ([yshift = -\gap]P.east); 
\draw[Myarrow]([xshift  =  \gap]Q.south) --  ([xshift =  \gap]C.north);
\draw[Myarrow]([xshift  = -\gap]C.north) --   ([xshift = -\gap]Q.south);
\draw[Myarrow,dashed] ([xshift = -\gap]B.north) -- ([xshift =  -\gap]P.south);
\draw[Myarrow] ([xshift =  \gap]P.south) -- ([xshift =  \gap]B.north);
\draw[Myarrow] ([yshift =  \gap]A.east)   --   ([yshift =  \gap]B.west) ;
\draw[Myarrow] ([yshift = -\gap]B.west)  --   ([yshift = -\gap]A.east); 
\draw[Myarrow] ([yshift =  \gap]B.east)  --  ([yshift =  \gap]C.west) ;
\draw[Myarrow] ([yshift = -\gap]C.west) --   ([yshift = -\gap]B.east); 

\draw[Myarrow] ([xshift =  \gap]O.south) -- ([xshift =  \gap]A.north);

\draw[Myarrow]([xshift  =  \gap]C.south) --  ([xshift =  \gap]Z.north);
\draw[Myarrow]([xshift  = -\gap]Z.north) --   ([xshift = -\gap]C.south);
\draw[Myarrow,dashed] ([xshift = -\gap]Y.north) -- ([xshift =  -\gap]B.south);
\draw[Myarrow] ([xshift =  \gap]B.south) -- ([xshift =  \gap]Y.north);
\draw[Myarrow] ([yshift =  \gap]X.east)   --   ([yshift =  \gap]Y.west) ;
\draw[Myarrow] ([yshift = -\gap]Y.west)  --   ([yshift = -\gap]X.east); 
\draw[Myarrow] ([yshift =  \gap]Y.east)  --  ([yshift =  \gap]Z.west) ;
\draw[Myarrow] ([yshift = -\gap]Z.west) --   ([yshift = -\gap]Y.east); 
\end{tikzpicture}
\label{eq:gammaconvalphathetaLKCD1}
\end{equation}
In this case, $\alpha=1/2$ gives
\begin{equation}
\begin{tikzpicture}[auto,scale=2.0, baseline=(current  bounding  box.center)]
\newcommand*{\size}{\scriptsize}%
\newcommand*{\gap}{.2ex}%
\newcommand*{\width}{2.5}%
\newcommand*{\height}{1.0}%

\node (O) at (0,0) {$\dfrac{1}{\pi}\mathrm{Im}\left\{\widetilde h_{1/2,\theta}(e^{-i\pi}t|\mu)\right\}$};
\node (P) at ($(O)+(\width,0)$)  {$h_{1/2,\theta}(x|\mu)$};
\node (Q) at ($(P)+(\width,0)$)  {$\widetilde h_{1/2,\theta}(s|\mu)$};
\node (A) at ($(O)-(0,\height)$) {$\mu\frac{\sin\pi\theta}{\pi}\, t^{\theta-1}\star \beta_{1/2,1/2-\theta}(t)$};
\node (B) at ($(A)+(\width,0)$)  {$\mu\{g_{1/2,1}\star g_{1/2-\theta,0}\}(x)$};
\node (C) at ($(B)+(\width,0)$) {$\mu\dfrac{s^\theta}{\sqrt{s(1+s)}}$};
\node (X) at ($(A)-(0,\height)$) {$\left(\frac{1}{2} - \theta\right)\delta(t)+\frac{1}{2}\delta(t-1)$};
\node (Y) at ($(X)+(\width,0)$) {$\left(\frac{1}{2} -\theta\right)+\frac{1}{2}e^{-x}$};
\node (Z) at  ($(Y)+(\width,0)$) {$\dfrac{1/2 -\theta}{s}+\dfrac{1/2}{1+s}$};   

\draw[Myarrow] ([yshift =  \gap]O.east)   --   ([yshift =  \gap]P.west) ;
\draw[Myarrow] ([yshift = -\gap]P.west)  --   ([yshift = -\gap]O.east); 
\draw[Myarrow] ([yshift =  \gap]P.east)   --   ([yshift =  \gap]Q.west) ;
\draw[Myarrow] ([yshift = -\gap]Q.west)  --   ([yshift = -\gap]P.east); 
\draw[Myarrow]([xshift  =  \gap]Q.south) --  ([xshift =  \gap]C.north);
\draw[Myarrow]([xshift  = -\gap]C.north) --   ([xshift = -\gap]Q.south);
\draw[Myarrow,dashed] ([xshift = -\gap]B.north) -- ([xshift =  -\gap]P.south);
\draw[Myarrow] ([xshift =  \gap]P.south) -- ([xshift =  \gap]B.north);
\draw[Myarrow] ([yshift =  \gap]A.east)   --   ([yshift =  \gap]B.west) ;
\draw[Myarrow] ([yshift = -\gap]B.west)  --   ([yshift = -\gap]A.east); 
\draw[Myarrow] ([yshift =  \gap]B.east)  --  ([yshift =  \gap]C.west) ;
\draw[Myarrow] ([yshift = -\gap]C.west) --   ([yshift = -\gap]B.east); 

\draw[Myarrow] ([xshift =  \gap]O.south) -- ([xshift =  \gap]A.north);

\draw[Myarrow]([xshift  =  \gap]C.south) --  ([xshift =  \gap]Z.north);
\draw[Myarrow]([xshift  = -\gap]Z.north) --   ([xshift = -\gap]C.south);
\draw[Myarrow] ([yshift =  \gap]X.east)   --   ([yshift =  \gap]Y.west) ;
\draw[Myarrow] ([yshift = -\gap]Y.west)  --   ([yshift = -\gap]X.east); 
\draw[Myarrow] ([yshift =  \gap]Y.east)  --  ([yshift =  \gap]Z.west) ;
\draw[Myarrow] ([yshift = -\gap]Z.west) --   ([yshift = -\gap]Y.east); 
\end{tikzpicture}
\label{eq:gammaconvalfahalfthetaLKCD}
\end{equation}

In the latter discussion
\begin{align}
-\frac{\widetilde h_{\alpha,\theta}^\prime(s|\mu)}{\widetilde h_{\alpha,\theta}(s|\mu)}
  &= \mu\dfrac{(1+s)^{\alpha-1}}{s^{\alpha-\theta}}
  \equiv \mu\,\rho_{\alpha,\theta}(s)  \\
    \implies\quad  \widetilde h_{\alpha,\theta}(s|\mu)
  &= \exp\left\{-\mu\int_0^s \rho_{\alpha,\theta}(t) dt \right\} 
\label{eq:h3param1}
 \end{align}

We note that  many densities that have arisen in our investigation of  the convolution of two gamma densities  
can be linked to a variety of other probabilistic  studies.
For example, 
$\beta_{1-\alpha,n\alpha+\theta}$ for $\theta>-\alpha$ and integer $n>0$ arises in 
the construction of the 
 Pitman-Yor  or two-parameter Poisson-Dirichlet  distribution ${\rm PD}(\alpha,\theta)$~\cite{PitmanYor},
which extends the original  one-parameter formulation ${\rm PD}(\alpha)$ 
 due to Kingman~\cite{Kingman,KingmanBook}.
 
 
 It is  only natural to explore  the generation of further densities through higher  convolutions of known  densities,
 or sums of such densities treated as Levy densities of some higher ID densities.
 We defer such further investigation to a separate study.
 Instead, we turn next to the study of  mixtures of stable densities and explore associated LKCD representations.

\section{Mixtures of Stable Densities}
\label{sec:stablemixtures}

As discussed above, the convolution of ID densities is also ID. 
It turns out that the sum of ID densities can also be ID, although this is by no means obvious from the ID representation by itself.
An example is a mixture (weighted sum or integral) of exponentials, which was shown to be ID by Steutel~\cite{Steutel1}.
We now broaden the foregoing discussion by allowing the scale parameter $y>0$  of the stable density $f_\alpha(x|y)$
($0<\alpha<1$) to be  a variable governed by a distribution with density $f(y|\mu)$.
In principle, $f$ might be any density, 
possibly involving several parameters. 
We  retain explicit dependence on at least one parameter  $\mu$ in anticipation of choosing an infinitely divisible $f(y|\mu)$ 
of the generic form  illustrated in the LKCD of Figure~\ref{fig:LKCD}.

The two-dimensional  joint density  of $x$ and $y$ is 
$\Pr(x,y|\mu)=\Pr(x|y)\Pr(y|\mu)= f_\alpha(x|y)f(y|\mu)$.
Hence the one-dimensional  marginal density of $x$ is 
\begin{align}
m_\alpha(x|\mu)\equiv\Pr(x|\mu) &= \int_0^\infty \Pr(x,y|\mu) dy = \int_0^\infty f_\alpha(x|y) f(y|\mu) dy
\label{eq:stablemixture}
\end{align}
This may be regarded  as  a  weighted mixture of stable densities  at different scales,
with mixing density $f(y|\mu)$.
Since $f_\alpha(x|y)$ 
has   $\widetilde f_\alpha(s|y)=\exp(-ys^\alpha)$, $m_\alpha(x|\mu)$ has Laplace  transform
\begin{align}
\widetilde m_\alpha(s|\mu) & = \int_0^\infty \widetilde f_\alpha(s|y) f(y|\mu) dy
    = \int_0^\infty e^{-y s^\alpha} f(y|\mu) dy
    = \widetilde f(s^\alpha|\mu)
\label{eq:stablemixturelaplace}
\end{align}
where $\widetilde f(s|\mu)$ is the Laplace  transform of $f(y|\mu)$.
Although $0<\alpha<1$ for a  stable distribution on a positive variable, we can accommodate $\alpha=1$ 
by  defining $f_{\alpha=1}(x|y)=\delta(x-y)$ with 
Laplace  transform $\exp(-ys)$, so that $m_{\alpha=1}$ reproduces $f$,
{\it i.e.}\ $m_{\alpha=1}(x|\mu)=f(x|\mu)$ and  $\widetilde m_{\alpha=1}(s|\mu) = \widetilde f(s|\mu)$.

An alternative approach that leads to~(\ref{eq:stablemixturelaplace}) is to consider the density of the product $Y^{1/\alpha}X$
where $X$ and $Y$ are independent variables with densities $f_\alpha(x)$ and $f(y)$ respectively, as discussed by 
Feller~\cite{Feller2}  (p463, Problem~10) and 
Bondesson~\cite{Bondesson} (p38, Example~3.2.4 for the case where $f$ is the  gamma density).
The study of distributions of products of independent variables is  a recurring theme in 
James~\cite{MR2676940} and James {\it et al.}~\cite{JamesRoynetteYor}.

Since  $f_\alpha$ is not available in closed form for general $0<\alpha<1$, neither is $m_\alpha$. 
But 
$\widetilde m_\alpha$ depends solely on the availability  of  $\widetilde f$.
If $f$ is ID then so is $m_\alpha$, as discussed under subordination and completely monotone functions in Feller~\cite{Feller2}~p451. 
Furthermore, Bondesson~\cite{Bondesson} (p41, Theorem~3.3.2) proved that if $f$ is GGC  then so is $m_\alpha$,
{\it i.e.}\ if $\widetilde f(s|\mu)$ is the Laplace  transform of a GCC,  so is  $\widetilde f(s^\alpha|\mu)$.
 In this case,  the availability  of  $\widetilde f(s^\alpha|\mu)$
 enables the generation of   Thorin densities that may not otherwise be readily identifiable as legitimate densities.
 Henceforth we shall confine interest to the GGC case.
 
 We always take $f_\alpha(x|\mu)$ ($0<\alpha<1$)  to denote the stable density.
 If we also reserve $\rho_\alpha(x)$ to be the  
 $\rho$-density of $f_\alpha(x|\mu)$, we need a
 different symbol, $r_\alpha$ say, for the $\rho$-density of $m_\alpha(x|\mu)$ -- {\it i.e.}\
 $r_\alpha(x)$ is to $m_\alpha(x|\mu)$ what $\rho(x)$
is to the general ID density $f(x|\mu)$ 
(and what  $\rho_\alpha(x)$  is to  $f_\alpha(x|\mu)$). 
Hence $\{\widetilde\rho(s), \widetilde r_\alpha(s)\}$, the Laplace transforms  of  $\{\rho(x), r_\alpha(x)\}$ respectively, are given by
 \begin{align}
\mu \widetilde\rho(s) &= -\frac{\widetilde f\,^\prime(s|\mu)}{\widetilde f(s|\mu)} 
\quad{\rm and}\quad 
\mu\,\widetilde r_\alpha(s) = -\frac{\widetilde m_\alpha\,^\prime(s|\mu)}{\widetilde m_\alpha(s|\mu)} 
  = -\frac{\widetilde f\,^\prime(s^\alpha|\mu)}{\widetilde f(s^\alpha|\mu)} 
\end{align}
To be explicit 
 \begin{align}
\mu\widetilde r_\alpha(s) 
&= -\frac{1}{\widetilde f(s^\alpha|\mu)}\frac{d}{ds}{\widetilde f(s^\alpha|\mu)}
   = -\frac{\alpha s^{\alpha-1}}{\widetilde f(s^\alpha|\mu)}\frac{d}{d s^\alpha}{\widetilde f(s^\alpha|\mu)}
   =\mu \alpha s^{\alpha-1}  \widetilde\rho(s^\alpha)
\end{align}
The limit forms for $\{\rho(x), r_\alpha(x)\}$ are
\begin{align}
 \mu\,\rho(x) &= \lim_{n\to\infty}n\,x\,f(x|\tfrac{\mu}{n})  \\
{\rm and}\quad \mu\, r_\alpha(x) 
&= \lim_{n\to\infty}n\,x\,m_\alpha(x|\tfrac{\mu}{n}) \\
  &= \phantom{\mu\,} x \int_0^\infty f_\alpha(x|y) \lim_{n\to\infty}n\,f(y|\tfrac{\mu}{n}) dy \\
  &= \mu\, x \int_0^\infty f_\alpha(x|y)\,y^{-1} \rho(y) dy
\label{eq:stablemixturerho} 
\end{align}

 Accordingly, for GGC $f(x|\mu)$,  $m_\alpha(x|\mu)$  has the following  GGC LKCD:
\begin{equation}
\begin{tikzpicture}[auto,scale=2.0, baseline=(current  bounding  box.center)]
\newcommand*{\size}{\scriptsize}%
\newcommand*{\gap}{.2ex}%
\newcommand*{\width}{2.25}%
\newcommand*{\height}{1.0}%

\node (O) at (0,0)  {$\dfrac{1}{\pi}\,{\rm Im}\left\{\widetilde f(e^{-i\pi\alpha}\,t^\alpha|\mu)\right\}$};
\node (P) at ($(O)+(\width,0)$)   {$m_\alpha(x|\mu)$};
\node (Q) at ($(P)+(\width,0)$)   {$\widetilde f(s^\alpha|\mu)$};
\node (A) at ($(O)-(0,\height)$) {$\dfrac{\mu\alpha}{\pi} {\rm Im}\left\{(e^{-i\pi}t)^{\alpha-1}\widetilde \rho(e^{-i\pi\alpha}\,t^\alpha)\right\}$};
\node (B) at ($(A)+(\width,0)$) {$\mu\, r_\alpha(x)$};
\node (C) at ($(B)+(\width,0)$) {$\mu\alpha s^{\alpha-1}\,\widetilde\rho(s^\alpha)$};    

\draw[Myarrow] ([yshift =  \gap]O.east)   --   ([yshift =  \gap]P.west) ;
\draw[Myarrow] ([yshift = -\gap]P.west)  --   ([yshift = -\gap]O.east); 
\draw[Myarrow] ([yshift =  \gap]P.east)   --   ([yshift =  \gap]Q.west) ;
\draw[Myarrow] ([yshift = -\gap]Q.west)  --   ([yshift = -\gap]P.east); 
\draw[Myarrow]([xshift  =  \gap]Q.south) --  ([xshift =  \gap]C.north);
\draw[Myarrow]([xshift  = -\gap]C.north) --   ([xshift = -\gap]Q.south);
\draw[Myarrow,dashed] ([xshift = -\gap]B.north) -- ([xshift =  -\gap]P.south);
\draw[Myarrow] ([xshift =  \gap]P.south) -- ([xshift =  \gap]B.north);
\draw[Myarrow] ([yshift =  \gap]A.east)   --   ([yshift =  \gap]B.west) ;
\draw[Myarrow] ([yshift = -\gap]B.west)  --   ([yshift = -\gap]A.east); 
\draw[Myarrow] ([yshift =  \gap]B.east)  --  ([yshift =  \gap]C.west) ;
\draw[Myarrow] ([yshift = -\gap]C.west) --   ([yshift = -\gap]B.east); 
\draw[Myarrow]([xshift  =  \gap]O.south) --  ([xshift =  \gap]A.north);
\end{tikzpicture}
\label{eq:hstableLKCD}
\end{equation}
The assertion that the bottom left node of~(\ref{eq:hstableLKCD}) is a density, 
even though the positivity of the expression  may not be obvious from mere inspection,
 restates  Bondesson's Theorem~3.3.2, with an  overlay of the GGC LKCD theme of this paper.

 \section{Stable Mixing Density} 
 \label{sec:stable}
 The stable-stable mixture density is discussed by Feller at various places in~\cite{Feller2} (pp176,~348,~452).
 Let $f(y|\mu)=f_\beta(y|\mu)  \implies \widetilde f(s|\mu)=\exp(-\mu s^\beta)$.
 We relabel $m_\alpha$ as $ m_{\alpha,\beta}$.
 Then $\widetilde  m_{\alpha,\beta}(s|\mu)= \widetilde f_\beta(s^\alpha|\mu)=\exp(-\mu s^{\alpha\beta})$,
  as can readily be verified.
  Hence $m_{\alpha,\beta} \equiv  f_{\alpha\beta}$, the stable density with parameter $\alpha\beta$.
  The GGC LKCD of $f_{\alpha\beta}(x|\mu)$ is 
  given in~(\ref{eq:extended3level3columnstableLKCD}), with $\alpha$ replaced by $\alpha\beta$.

The integral representations~(\ref{eq:stablemixture}) and~(\ref{eq:stablemixturerho}) take the form
  \begin{align}
      m_{\alpha,\beta}(x|\mu) \equiv f_{\alpha\beta}(x|\mu)  &= \int_0^\infty f_\alpha(x|y) f_\beta(y|\mu) dy
  \label{eq:stablestablemixture} \\
   r_{\alpha,\beta}(x) \equiv \rho_{\alpha\beta}(x) 
     &= x \int_0^\infty f_\alpha(x|y) y^{-1} \rho_\beta(y) dy
  \label{eq:stablestablerho} 
  \end{align}
 
For $\alpha=\beta=1/2$, the known density $f_{1/2}$ of~(\ref{eq:3levelstableLKCD}) induces an integral representation for $f_{1/4}$
  \begin{align}
     f_{1/4}(x|\mu)    &= \int_0^\infty f_{1/2}(x|y) f_{1/2}(y|\mu) dy \\
     &=  \int_0^\infty  \frac{y}{2\sqrt{\pi x^3}}\, e^{-y^2/4x}\,\frac{\mu}{2\sqrt{\pi y^3}}\, e^{-\mu^2/4y} \, dy \nonumber \\
     &=   \frac{\mu}{4\pi\sqrt{x^3}}\int_0^\infty  \frac{1}{\sqrt{y}}\, e^{-\mu^2/4y-y^2/4x} \, dy
  \label{eq:stablequarter} 
  \end{align}
  Berberan-Santos~\cite{Berberan-Santos} derived (\ref{eq:stablequarter})  for $\mu=1$ through Laplace inversion. 

 The corresponding $\rho_{\alpha\beta}(x)$, with  $\alpha=\beta=1/2$,  is
 \begin{align}
 \frac{\alpha\beta x^{-\alpha\beta}}{\Gamma(1-\alpha\beta)} &= \frac{1}{4}\, \frac{x^{-1/4}}{\Gamma(\frac{3}{4})}
\intertext{  which is reproduced, as it should be, by the limit}
 \lim_{n\to\infty} n x  f_{1/4}(x|\tfrac{\mu}{n})  &= \frac{\mu}{4\pi\sqrt{x}} \int_0^\infty y^{-1/2}\, e^{-y^2/4x} \, dy \\
    &=  \frac{\mu}{8\pi\sqrt{x}} \int_0^\infty y^{-3/4} e^{-y/4x} \, dy \quad (y^2\to y) \\
    &= \frac{\mu}{8\pi\sqrt{x}}\, \Gamma(\tfrac{1}{4}) \,(4x)^{1/4} \\
    &= \frac{\mu}{4}\, \frac{x^{-1/4}}{\Gamma(\frac{3}{4})}
  \end{align}
since $\Gamma(\tfrac{1}{4})\Gamma(\tfrac{3}{4}) = \pi\sqrt{2}$ by the Legendre duplication  formula. 

As noted earlier, $f_{1/4}$ is of particular interest in physics.
We may then use $\alpha=1/2$ and $\beta=1/4$ (or the other way round) 
in~(\ref{eq:stablestablemixture}) to generate $f_{1/8}$, although the integral  representation will inevitably be more complex.

 \section{Gamma Mixing Density}
 \label{sec:hgamma}
 Let $f(y|\mu)$ be the gamma density $g_{\mu,\lambda}(y)$ of (\ref{eq:gamma}).
 Then 
 $m_{\alpha,\lambda}(x|\mu)$ has Laplace transform 
 \begin{align}
\widetilde m_{\alpha,\lambda}(s|\mu) &=  \widetilde g_{\mu,\lambda}(s^\alpha) = \frac{\lambda^\mu}{(\lambda+s^\alpha)^\mu}
\label{eq:fracgamma}
 \end{align}  
$\alpha=1$ reproduces the Laplace transform of the gamma density $\widetilde g_{\mu,\lambda}(s)$.
 For $0<\alpha<1$~(\ref{eq:fracgamma}) is the Laplace transform of what is  known as 
 the  fractional gamma distribution, 
{\it e.g.}\ Di~Nardo {\it et al.}~\cite{DiNardo2021} (Definition~2.1).
It is also called the positive Linnik distribution.
It also seems reasonable to refer to it as the stable-gamma mixture distribution.

 The GGC LKCD of the fractional gamma density  $m_{\alpha,\lambda}(x|\mu)$ is
 \begin{equation}
\begin{tikzpicture}[auto,scale=2.0, baseline=(current  bounding  box.center)]
\newcommand*{\size}{\scriptsize}%
\newcommand*{\gap}{.2ex}%
\newcommand*{\width}{2.25}%
\newcommand*{\height}{1.25}%

\node (O) at (0,0)  {$\dfrac{1}{\pi}\,{\rm Im}\left\{\dfrac{\lambda^\mu}{(\lambda+e^{-i\pi\alpha}\,t^\alpha)^\mu}\right\}$};
\node (P) at ($(O)+(\width,0)$)   {$m_{\alpha,\lambda}(x|\mu)$};
\node (Q) at ($(P)+(\width,0)$)   {$\dfrac{\lambda^\mu}{(\lambda+s^\alpha)^\mu}$};
\node (A) at ($(O)-(0,\height)$) {$\dfrac{\mu\alpha}{\pi} {\rm Im}\left\{\dfrac{(e^{-i\pi}t)^{\alpha-1}}{\lambda+e^{-i\pi\alpha}\,t^\alpha}\right\}$};
\node (B) at ($(A)+(\width,0)$) {$\mu\, r_{\alpha,\lambda}(x)$};
\node (C) at ($(B)+(\width,0)$) {$\mu\alpha \dfrac{s^{\alpha-1}}{\lambda+s^\alpha}$};    
\draw[Myarrow] ([yshift =  \gap]O.east)   --   ([yshift =  \gap]P.west) ;
\draw[Myarrow] ([yshift = -\gap]P.west)  --   ([yshift = -\gap]O.east); 
\draw[Myarrow] ([yshift =  \gap]P.east)   --   ([yshift =  \gap]Q.west) ;
\draw[Myarrow] ([yshift = -\gap]Q.west)  --   ([yshift = -\gap]P.east); 
\draw[Myarrow]([xshift  =  \gap]Q.south) --  ([xshift =  \gap]C.north);
\draw[Myarrow]([xshift  = -\gap]C.north) --   ([xshift = -\gap]Q.south);
\draw[Myarrow,dashed] ([xshift = -\gap]B.north) -- ([xshift =  -\gap]P.south);
\draw[Myarrow] ([xshift =  \gap]P.south) -- ([xshift =  \gap]B.north);
\draw[Myarrow] ([yshift =  \gap]A.east)   --   ([yshift =  \gap]B.west) ;
\draw[Myarrow] ([yshift = -\gap]B.west)  --   ([yshift = -\gap]A.east); 
\draw[Myarrow] ([yshift =  \gap]B.east)  --  ([yshift =  \gap]C.west) ;
\draw[Myarrow] ([yshift = -\gap]C.west) --   ([yshift = -\gap]B.east); 
\draw[Myarrow]([xshift  =  \gap]O.south) --  ([xshift =  \gap]A.north);
\end{tikzpicture}
\label{eq:paraboliccylinderLKCD}
\end{equation}
The Thorin density in the bottom left node takes the explicit form
\begin{align}
\frac{\mu\alpha}{\pi} {\rm Im}\left\{\frac{(e^{-i\pi}t)^{\alpha-1}}{\lambda+e^{-i\pi\alpha}\,t^\alpha}\right\}
&= \frac{\mu\alpha}{\pi}  \frac{\lambda\,t^{\alpha-1}\sin\pi\alpha}{\lambda^2+2\lambda\,t^\alpha\cos\pi\alpha+t^{2\alpha}}
\label{eq:thorin1}
\end{align}
This was derived for $\lambda=1$ by Bondesson~\cite{Bondesson} (p38), as a consequence of Theorem~3.3.2~(p41).
As previously discussed, the theorem gives assurance that~(\ref{eq:thorin1}) is  a valid density.

 What form does $m_{\alpha,\lambda}(x|\mu)$ itself take?

\subsection{Geometric Series  Representation} 
\label{sec:geometricseries}
Let $m_{\alpha,\lambda}(x)\equiv m_{\alpha,\lambda}(x|\mu=1)$ with Laplace transform $\widetilde m_{\alpha,\lambda}(s)$,  
which we expand  as a geometric series 
 \begin{align}
\widetilde m_{\alpha,\lambda}(s)  &= \frac{\lambda}{\lambda+s^\alpha} 
      = \frac{\lambda/s^\alpha}{1+\lambda/s^\alpha}
      = -\sum_{k=1}^\infty  \left(-\frac{\lambda}{s^\alpha}\right)^{k}  \\
\implies 
m_{\alpha,\lambda}(x) &= -\sum_{k=1}^\infty (-\lambda)^k  \frac{x^{\alpha k-1}}{\Gamma(\alpha k)}
\label{eq:MittagLeflerdensity}      
 \end{align}  
 Let $F_{\alpha,\lambda}(x)$ be the distribution  with density $m_{\alpha,\lambda}(x)$, so that
 \begin{align}
F_{\alpha,\lambda}(x) &= \int_0^x m_{\alpha,\lambda}(y)dy \\
  &= -\sum_{k=1}^\infty \frac{(-\lambda\, x^\alpha)^k}{\Gamma(\alpha k+1)}
  = 1-E_\alpha(-\lambda\, x^\alpha)
\label{eq:MittagLeflerdistribution}      
\intertext{where $E_\alpha(x)$ is the Mittag-Leffler function}
    E_\alpha(x)  &=  \sum_{k=0}^\infty  \frac{x^{k}}{\Gamma(\alpha k+1)}  
\label{eq:MittagLeflerfunction}      
\intertext{Hence the density $m_{\alpha,\lambda}(x)\equiv m_{\alpha,\lambda}(x|\mu=1)$ can be written as}
m_{\alpha,\lambda}(x) &= F_{\alpha,\lambda}^{\,\prime}(x) = -E_\alpha^{\,\prime}(-\lambda\, x^\alpha)
\label{eq:MittagLeflerdensity1}      
\end{align}  

Feller~(\cite{Feller2}~p453) discussed the Mittag-Leffler function  in the context of Laplace transforms in two dimensions.
Pillai~\cite{Pillai} defined the Mittag-Leffler distribution as  
the  case  $F_{\alpha,\lambda=1}(x)$  (in the notation of this paper).
Hauboldt~{\it et al.}~\cite{Haubold2011MittagLefflerFA} gave a review of the Mittag-Leffler function. 
The general case  $F_{\alpha,\lambda}(x|\mu)$ with density $m_{\alpha,\lambda}(x|\mu)$
may be expressed in terms of the generalised (three parameter) Mittag-Leffler function,
also known as the Prabhakar function, as reviewed by Garra and Garrappa~\cite{Garra}.

In light of the foregoing discussion, we may plausibly  refer to $m_{\alpha,\lambda}(x)$ as the Mittag-Leffler density.
Commonly used though the Mittag-Leffler geometric  series representation may be, it is not a unique 
representation of the fractional gamma distribution.
We demonstrate an alternative perspective next that directly adheres to  the integral representation of   $m_{\alpha,\lambda}(x|\mu)$.

\subsection{Integral Representation of Fractional Gamma Density}
 \label{sec:fracgammaintrep}

The  integral representations~(\ref{eq:stablemixture}),~(\ref{eq:stablemixturerho}) are
  \begin{align}
     m_{\alpha,\lambda}(x|\mu) 
     &=  \frac{\lambda^\mu}{\Gamma(\mu)}  \int_0^\infty   f_\alpha(x|y)\, y^{\mu-1}\, e^{-\lambda y}  \, dy
  \label{eq:stablegammamixture}  \\
    r_{\alpha,\lambda}(x) 
     &=  x \int_0^\infty   f_\alpha(x|y)\, y^{-1}\, e^{-\lambda y}  \, dy
  \label{eq:stablegammamixturethorin} 
 \end{align}
Choosing $\alpha=1/2$ gives
\begin{align}
     m_{1/2,\lambda}(x|\mu)   
     &= \frac{\lambda^\mu}{2\Gamma(\mu)\sqrt{\pi x^3}}  \int_0^\infty   y^\mu \, e^{-y^2/4x-\lambda y} \, dy \\
     &= \sqrt{\frac{2}{\pi}}\, \mu\lambda^\mu\, (2x)^{\mu/2-1} \, e^{\lambda^2 x/2}\,  D_{-\mu-1}(\lambda\sqrt{2x}\,)
  \label{eq:stablehalfgammamixture} 
  \end{align}
  where $D_\mu(x)$ is the parabolic cylinder function~\cite{Gradshteyn} (p365,~3.462.1).
$D_\mu(x)$  arises in the solution of Laplace's equation by separation of variables  in parabolic cylinder coordinates. 
 Also
\begin{align}
r_{1/2,\lambda}(x) 
&= \frac{1}{2}\, \sqrt{\frac{2}{\pi}}\, e^{\lambda^2 x/2}\, D_{-1}(\lambda\sqrt{2x}\,)
     = \frac{1}{2}\, e^{\lambda^2x}\erfc(\lambda\sqrt{x})
  \label{eq:stablehalfgammamixturethorin} 
\end{align}
 The rightmost form follows from~\cite{Gradshteyn}~(p1030,~9.254.1), $\erfc(x)$ being the complementary error function.
 It is compatible with the Laplace transform route, with the aid of the Laplace transform
 pairs~\cite{AbrSteg}~(p1028, 29.3.114 and p1027, 29.3.90).
 The GGC LKCD of the fractional gamma density  $m_{1/2,\lambda}(x|\mu)$ is
 \begin{equation}
\begin{tikzpicture}[auto,scale=2.0, baseline=(current  bounding  box.center)]
\newcommand*{\size}{\scriptsize}%
\newcommand*{\gap}{.2ex}%
\newcommand*{\width}{2.25}%
\newcommand*{\height}{1.25}%

\node (O) at (0,0)  {$\dfrac{1}{\pi}\dfrac{{\rm Im}\left(\lambda+i\sqrt{t}\,\right)^\mu}{(\lambda^2+t)^\mu}$};
\node (P) at ($(O)+(\width,0)$)   {$m_{1/2,\lambda}(x|\mu)$};
\node (Q) at ($(P)+(\width,0)$)   {$\dfrac{\lambda^\mu}{(\lambda+\sqrt{s}\,)^\mu}$};
\node (A) at ($(O)-(0,\height)$) {$\dfrac{\mu}{2\pi}\dfrac{\lambda}{\sqrt{t}\left(\lambda^2+t\right)}$};
\node (B) at ($(A)+(\width,0)$) {$\dfrac{\mu}{2} e^{\lambda^2x}\erfc(\lambda\sqrt{x})$};
\node (C) at ($(B)+(\width,0)$) {$\dfrac{\mu}{2}\dfrac{1}{\sqrt{s}\left(\lambda+\sqrt{s}\,\right)}$};    
\draw[Myarrow] ([yshift =  \gap]O.east)   --   ([yshift =  \gap]P.west) ;
\draw[Myarrow] ([yshift = -\gap]P.west)  --   ([yshift = -\gap]O.east); 
\draw[Myarrow] ([yshift =  \gap]P.east)   --   ([yshift =  \gap]Q.west) ;
\draw[Myarrow] ([yshift = -\gap]Q.west)  --   ([yshift = -\gap]P.east); 
\draw[Myarrow]([xshift  =  \gap]Q.south) --  ([xshift =  \gap]C.north);
\draw[Myarrow]([xshift  = -\gap]C.north) --   ([xshift = -\gap]Q.south);
\draw[Myarrow,dashed] ([xshift = -\gap]B.north) -- ([xshift =  -\gap]P.south);
\draw[Myarrow] ([xshift =  \gap]P.south) -- ([xshift =  \gap]B.north);
\draw[Myarrow] ([yshift =  \gap]A.east)   --   ([yshift =  \gap]B.west) ;
\draw[Myarrow] ([yshift = -\gap]B.west)  --   ([yshift = -\gap]A.east); 
\draw[Myarrow] ([yshift =  \gap]B.east)  --  ([yshift =  \gap]C.west) ;
\draw[Myarrow] ([yshift = -\gap]C.west) --   ([yshift = -\gap]B.east); 
\draw[Myarrow]([xshift  =  \gap]O.south) --  ([xshift =  \gap]A.north);
\end{tikzpicture}
\label{eq:paraboliccylinderLKCD1}
\end{equation}

With the aid of~(\ref{eq:stablestablemixture}) and~(\ref{eq:stablestablerho}) in the stable-stable case above  we get
  \begin{align}
     m_{\alpha\beta,\lambda}(x|\mu) 
     &=  \frac{\lambda^\mu}{\Gamma(\mu)}  \int_0^\infty   f_{\alpha\beta}(x|y)\, y^{\mu-1}\, e^{-\lambda y}  \, dy \\
     &= \int_0^\infty  f_{\alpha}(x|u)
            \left\{ \frac{\lambda^\mu}{\Gamma(\mu)} \int_0^\infty   f_{\beta}(u|y)\, y^{\mu-1}\, e^{-\lambda y}  \, dy \right\} du \\
     &=  \int_0^\infty  f_{\alpha}(x|u)\, m_{\beta,\lambda}(u|\mu) du
\label{eq:malphabeta}    \\
  r_{\alpha\beta,\lambda}(x) 
       &=  x \int_0^\infty  f_{\alpha}(x|u)\, u^{-1} \, r_{\beta,\lambda}(u) du
\label{eq:ralphabeta}     
\end{align}
Hence we may use  the case $\alpha=1/2$ to induce the $\alpha=1/4$ case
  \begin{align}
     m_{1/4,\lambda}(x|\mu)  
     &= \int_0^\infty   f_{1/2}(x|z)\, m_{1/2,\lambda}(z|\mu) \, dz \\
     &=   \frac{1}{2\sqrt{\pi x^3}}\int_0^\infty  z\, e^{-z^2/4x} \, m_{1/2,\lambda}(z|\mu) \, dz
  \label{eq:stablequartergammamixture} \\
  r_{1/4,\lambda}(x)
   &=  \frac{1}{2\sqrt{\pi x}}\int_0^\infty   e^{-z^2/4x}\, r_{1/2,\lambda}(z) \, dz \\
   &=  \frac{1}{4\sqrt{\pi x}}\int_0^\infty   e^{-z^2/4x+\lambda^2z}\erfc(\lambda\sqrt{z}\,) \, dz 
  \end{align}
  
 The GGC LKCD of the fractional gamma density  $m_{1/4,\lambda}(x|\mu)$ is shown below
 \begin{equation}
\begin{tikzpicture}[auto,scale=2.0, baseline=(current  bounding  box.center)]
\newcommand*{\size}{\scriptsize}%
\newcommand*{\gap}{.2ex}%
\newcommand*{\width}{2.4}%
\newcommand*{\height}{1.25}%

\node (O) at (0,0)  {$\dfrac{1}{\pi}\,{\rm Im}\left\{\dfrac{\lambda^\mu}{(\lambda+e^{-i\pi/4}\,t^{1/4})^\mu}\right\}$};
\node (P) at ($(O)+(\width,0)$)   {$m_{1/4,\lambda}(x|\mu)$};
\node (Q) at ($(P)+(\width,0)$)   {$\dfrac{\lambda^\mu}{\left(\lambda+s^{1/4}\right)^\mu}$};
\node (A) at ($(O)-(0,\height)$) {$\dfrac{\mu}{4\sqrt{2}\pi}\dfrac{\lambda\,t^{-3/4}}{\left(\lambda^2+\sqrt{2}\,\lambda\,t^{1/4}+\sqrt{t}\,\right)}$};
\node (B) at ($(A)+(\width,0)$) {$\mu\, r_{1/4,\lambda}(x)$};
\node (C) at ($(B)+(\width,0)$) {$\dfrac{\mu}{4}\dfrac{1}{s^{3/4}\left(\lambda+s^{1/4}\right)}$};    
\draw[Myarrow] ([yshift =  \gap]O.east)   --   ([yshift =  \gap]P.west) ;
\draw[Myarrow] ([yshift = -\gap]P.west)  --   ([yshift = -\gap]O.east); 
\draw[Myarrow] ([yshift =  \gap]P.east)   --   ([yshift =  \gap]Q.west) ;
\draw[Myarrow] ([yshift = -\gap]Q.west)  --   ([yshift = -\gap]P.east); 
\draw[Myarrow]([xshift  =  \gap]Q.south) --  ([xshift =  \gap]C.north);
\draw[Myarrow]([xshift  = -\gap]C.north) --   ([xshift = -\gap]Q.south);
\draw[Myarrow,dashed] ([xshift = -\gap]B.north) -- ([xshift =  -\gap]P.south);
\draw[Myarrow] ([xshift =  \gap]P.south) -- ([xshift =  \gap]B.north);
\draw[Myarrow] ([yshift =  \gap]A.east)   --   ([yshift =  \gap]B.west) ;
\draw[Myarrow] ([yshift = -\gap]B.west)  --   ([yshift = -\gap]A.east); 
\draw[Myarrow] ([yshift =  \gap]B.east)  --  ([yshift =  \gap]C.west) ;
\draw[Myarrow] ([yshift = -\gap]C.west) --   ([yshift = -\gap]B.east); 
\draw[Myarrow]([xshift  =  \gap]O.south) --  ([xshift =  \gap]A.north);
\end{tikzpicture}
\label{eq:paraboliccylinderLKCD2}
\end{equation}
We can proceed to generate an integral representation for $\alpha=1/8$ with the aid of $\alpha=1/2, 1/4$.
We can similarly start from $\alpha=1/3$ combined with $\alpha=1/2$ to generate the  sequence 
for $\alpha=\{1/3, 1/6, 1/9,\dots\}$. 

 For general $\alpha$, we may revert to the  Pollard infinite series representation~(\ref{eq:pollardsum2})
 for $f_\alpha$ in the integral representations~(\ref{eq:stablegammamixture}) and (\ref{eq:stablegammamixturethorin}):
 \begin{align}
     m_{\alpha,\lambda}(x|\mu) 
    &=  -\frac{\lambda^\mu}{\pi\Gamma(\mu)} \sum_{k=0}^\infty \frac{(-1)^k}{k!} \sin(\pi k \alpha)\frac{\Gamma(k \alpha+1)}{x^{k \alpha+1}} 
              \int_0^\infty y^{\mu+k-1}  e^{-\lambda y}  \, dy \\
    &=  -\frac{1}{\pi\Gamma(\mu)} \sum_{k=0}^\infty \frac{(-1)^k}{k!} \sin(\pi k \alpha)
            \frac{\Gamma(k \alpha+1)\Gamma(\mu+k)}{\lambda^k x^{k \alpha+1}} 
\label{eq:pollardgammamixture} \\
{\rm and}\quad    r_{\alpha,\lambda}(x) 
    &=  - \frac{1}{\pi}\sum_{k=1}^\infty  \frac{(-1)^k}{k} \sin(\pi k \alpha)\frac{\Gamma(k \alpha+1)}{\lambda^k x^{k \alpha}} 
    =  -\alpha \sum_{k=1}^\infty  \frac{(-\lambda x^\alpha)^{-k}}{\Gamma(1-k \alpha)} 
\label{eq:ML1}  
\end{align}

The geometric series Mittag-Leffler representation  is, of course, fully consistent with the integral representation, including the 
parabolic cylinder function representation. 
Explicitly, for $\alpha=1/2$
\begin{align}
2 \, r_{1/2,\lambda}(x &)=\sqrt{\frac{2}{\pi}}\,e^{\lambda^2 x/2}\,D_{-1}(\lambda\sqrt{2x})=E_{1/2}(-\lambda\sqrt{x})
=e^{\lambda^2 x}\,\erfc(\lambda \sqrt{x}) \\
m_{1/2,\lambda}(x) &=   \frac{\lambda}{\sqrt{\pi x}}  \, e^{\lambda^2 x/2}\,  D_{-2}(\lambda\sqrt{2x}\,) \\
   &= -\frac{d}{dx}\left(e^{\lambda^2 x}\,\erfc(\lambda \sqrt{x})\right) 
   =  -E_{1/2}^{\,\prime}(-\lambda\sqrt{x}) 
\intertext{where, by~\cite{Gradshteyn} (p1030,~9.254.2)} 
D_{-2}(\lambda\sqrt{2x}) &= e^{-\lambda^2 x/2}- \lambda \sqrt{\pi x}\,  e^{\lambda^2x/2} \erfc(\lambda\sqrt{x})
\end{align}

\subsection{Convolution Revisited}
\label{sec:convmix}
Finally, we note that the fractional gamma mixture  may  be looked upon as arising from a  convolution of two densities:
$g_{1-\alpha,0}(x)$ and $m_{\alpha,\lambda}(x)\equiv m_{\alpha,\lambda}(x|\mu=1)$,
For brevity, we define 
\begin{align}
 \{g\star m\}_{\alpha,\lambda}(x) &\equiv g_{1-\alpha,0}(x)\star m_{\alpha,\lambda}(x)/\lambda \\
\mathrm{and}\quad 
\phi_{\alpha,\lambda}(t) &\equiv \frac{\alpha}{\pi} {\rm Im}\left\{\frac{(e^{-i\pi}t)^{\alpha-1}}{\lambda+e^{-i\pi\alpha}\,t^\alpha}\right\}
= \frac{\alpha}{\pi}  \frac{\lambda\,t^{\alpha-1}\sin\pi\alpha}{\lambda^2+2\lambda\,t^\alpha\cos\pi\alpha+t^{2\alpha}}
\end{align}
$\phi_{\alpha,\lambda}(t)$ is the Thorin density of~(\ref{eq:thorin1}).

Hence, in  keeping with the  convolution of two gamma densities  studied earlier,
we have the followed \textquote{convolution-centred} GGC LKCD for the fractional gamma density $m_{\alpha,\lambda}(x|\mu)$

\begin{equation}
\begin{tikzpicture}[auto,scale=2.0, baseline=(current  bounding  box.center)]
\newcommand*{\size}{\scriptsize}%
\newcommand*{\gap}{.2ex}%
\newcommand*{\width}{2.5}%
\newcommand*{\height}{1.0}%

\node (O) at (0,0)  {$\dfrac{1}{\pi}\,{\rm Im}\left\{\dfrac{\lambda^\mu}{(\lambda+e^{-i\pi\alpha}\,t^\alpha)^\mu}\right\}$};
\node (P) at ($(O)+(\width,0)$)   {$m_{\alpha,\lambda}(x|\mu)$};
\node (Q) at ($(P)+(\width,0)$)   {$\dfrac{\lambda^\mu}{(\lambda+s^\alpha)^\mu}$};
\node (A) at ($(O)-(0,\height)$) {$\mu\, \phi_{\alpha,\lambda}(t)$};
\node (B) at ($(A)+(\width,0)$) {$\mu\alpha \{g\star m\}_{\alpha,\lambda}(x)$};
\node (C) at ($(B)+(\width,0)$) {$\mu\alpha \dfrac{s^{\alpha-1}}{\lambda+s^\alpha}$};    

\draw[Myarrow] ([yshift =  \gap]O.east)   --   ([yshift =  \gap]P.west) ;
\draw[Myarrow] ([yshift = -\gap]P.west)  --   ([yshift = -\gap]O.east); 
\draw[Myarrow] ([yshift =  \gap]P.east)   --   ([yshift =  \gap]Q.west) ;
\draw[Myarrow] ([yshift = -\gap]Q.west)  --   ([yshift = -\gap]P.east); 
\draw[Myarrow]([xshift  =  \gap]Q.south) --  ([xshift =  \gap]C.north);
\draw[Myarrow]([xshift  = -\gap]C.north) --   ([xshift = -\gap]Q.south);
\draw[Myarrow,dashed] ([xshift = -\gap]B.north) -- ([xshift =  -\gap]P.south);
\draw[Myarrow] ([xshift =  \gap]P.south) -- ([xshift =  \gap]B.north);
\draw[Myarrow] ([yshift =  \gap]A.east)   --   ([yshift =  \gap]B.west) ;
\draw[Myarrow] ([yshift = -\gap]B.west)  --   ([yshift = -\gap]A.east); 
\draw[Myarrow] ([yshift =  \gap]B.east)  --  ([yshift =  \gap]C.west) ;
\draw[Myarrow] ([yshift = -\gap]C.west) --   ([yshift = -\gap]B.east); 
\draw[Myarrow]([xshift  =  \gap]O.south) --  ([xshift =  \gap]A.north);

\node (X) at ($(A)-(0,\height)$)  {$(1-\alpha)\delta(t)+\phi_{\alpha,\lambda}(t)$};
\node (Y) at ($(X)+(\width,0)$) {$1-\alpha+\alpha \{g\star m\}_{\alpha,\lambda}(x)$};
\node (Z) at  ($(Y)+(\width,0)$) {$\dfrac{1-\alpha}{s}+\alpha \dfrac{s^{\alpha-1}}{\lambda+s^\alpha}$};   

\draw[Myarrow]([xshift  =  \gap]C.south) --  ([xshift =  \gap]Z.north);
\draw[Myarrow]([xshift  = -\gap]Z.north) --   ([xshift = -\gap]C.south);
\draw[Myarrow,dashed] ([xshift = -\gap]Y.north) -- ([xshift =  -\gap]B.south);
\draw[Myarrow] ([xshift =  \gap]B.south) -- ([xshift =  \gap]Y.north);
\draw[Myarrow] ([yshift =  \gap]X.east)   --   ([yshift =  \gap]Y.west) ;
\draw[Myarrow] ([yshift = -\gap]Y.west)  --   ([yshift = -\gap]X.east); 
\draw[Myarrow] ([yshift =  \gap]Y.east)  --  ([yshift =  \gap]Z.west) ;
\draw[Myarrow] ([yshift = -\gap]Z.west) --   ([yshift = -\gap]Y.east); 
\end{tikzpicture}
\label{eq:convstableLKCD}
\end{equation}

This demonstrates a unifying theme of the convolution of two densities for all objects studied in this paper, despite their apparent diversity.

\section{Discussion}
\label{sec:discussion}
It  is worth restating the objective of this paper, building upon the introductory remarks.
 The  novelty of the paper is  primarily one of perspective and representation rather than discovery of new, previously undocumented, ID distributions.
 Such novelty  lies 
 in representing  ID densities, known or novel,  as nodes of a 
 commutative diagram, with 
 \vspace{-10pt}
 \begin{myitemize}
 \item the Laplace  transform or its inverse as horizontal  connections of the diagram 
\item a  limiting process or logarithmic derivative as downward connections 
\item compound Poisson sum or integral of logarithmic derivative as upward connections.
 \end{myitemize}
  \vspace{-10pt}
Such visual representation has  vastly contributed to our own appreciation of the coherence and  connection amongst  densities 
that appear  disparate  at first.
 

 The typical goal of mathematical research is the quest for  generality. 
A case in point is the representation of the stable density in terms of  the very general Meijer G-function 
(Penson and G\'{o}rska~\cite{PensonGorska}).
By contrast, the approach of this paper has been one of conceptual simplicity, such as the convolution of two densities,
and then exploring the generality that may flow from such simplicity.

Accordingly, we have  introduced only one basic ID object, the gamma density $g_{\mu,\lambda}(x)$,
 which gave rise to the gamma LKCD~(\ref{eq:gammaLKCD}).
 Then, rather than introducing the stable  density as a new object, we inferred it from the gamma case by growing the gamma LKCD 
 upward from  $g_{1-\alpha,0}(x)=x^{-\alpha}/\Gamma(1-\alpha)$ ($0<\alpha<1$).
 We then also extended the LKCD to the left since $x^{-\alpha}$ is the Laplace transform of $g_{\alpha,0}(t)=t^{\alpha-1}/\Gamma(\alpha)$,
 to form the LKCD of a generalised gamma convolution, the GGC LKCD.
 We then turned to the convolution of two gamma densities, thereby entering the world of  beta and Bessel densities
 via the  confluent hypergeometric function.
 That all these objects  can be seen to arise from the gamma density or the   convolution of two gamma densities
 (which is not simply another gamma density if the decays are different), 
 can easily pass unappreciated, leading to a less joined-up conversation about them than might otherwise be the case.
 We take conceptual simplicity and the intrinsically joined-up commutative diagram to be  intimately  related fundamentals.

Another key perspective that we explored is the mixture of stable densities. 
If the mixing density is itself stable, the mixture is also  a stable density of index $\alpha\beta$ generated from stable densities
of separate indices $\alpha$ and $\beta$.
A gamma mixing density led to the fractional gamma density.
The latter  is routinely represented in terms of the Mittag-Leffler function.
We also demonstrated  an integral representation that gives an expression involving the parabolic cylinder function  
for the stable index $\alpha=1/2$.
To our awareness, the parabolic cylinder function representation of the fractional gamma distribution has not been reported in the literature.
As a solution of Laplace's equation in cylindrical coordinates, the parabolic cylinder function perspective  suggests a 
link between the fractional gamma distribution and the study of  harmonic functions in potential theory.

\section{Conclusion and Future Work}
\label{sec:conclusion}
We have introduced  a commutative diagram visualisation of  infinitely divisible (ID) distributions (or their densities, to be  precise).
We referred to this novel representation as the L\'{e}vy-Khintchine commutative diagram  (LKCD).

Much remains to be explored.
Notably, the ID densities studied here can form the basis for the construction of  multivariate distributions. 
The simplest case is the Dirichlet  distribution, which Ferguson~\cite{Ferguson1,Ferguson2} constructed as follows.
Let $\scrG(\mu,\lambda)$ denote the gamma distribution with shape $\mu$ and decay $\lambda$
and let $\{X_i\sim\scrG(\mu_i,\lambda)\}$ be $ N$  independently  distributed  variables  
with individual  gamma distributions that may have different shapes $\{\mu_i\}$ but share a common decay $\lambda$.
The distribution of the sum $X=\sum_{i=1}^N X_i$ is a convolution of the individual gamma distributions 
and thus also a gamma distribution whose shape is a sum of the individual shapes with the same  shared decay.
Then the multivariate  distribution on the normalised variables $\{X_i/X\}$ 
is known as  the Dirichlet  distribution. 
It is defined on the $(N-1)$-dimensional probability simplex  
in the $N$-dimensional space defined by the independent $\{X_i\}$.
Notably, it depends only on the $\{\mu_i\}$ and is independent  of the shared decay $\lambda$.


Alternatives  to Dirichlet are possible but more complex.
For instance, Sibisi and Skilling~\cite{SibisiSkilling} suggested the convolution of two gamma distribution 
$\{X_i\sim\scrG(\mu_i,\lambda_1)\star \scrG(\nu_i,\lambda_2)\}$,
which they referred to as the  supergamma distribution and
the induced  normalised distribution the superDirichlet distribution.
Di~Nardo {\it et al.}~\cite{DiNardo2021} explored  $\{X_i\sim\scrF\scrG(\mu_i,\alpha,\lambda)\}$ 
where $\scrF\scrG(\mu,\alpha,\lambda)$ is the  fractional gamma distribution, which they expressed in terms of the 
three parameter Mittag-Leffler (Prabhakar) function. 
 They referred to the associated  normalised distribution as the  fractional generalisation of the Dirichlet distribution.
 Favaro {\it et al.}~\cite{Favaro} discussed the general case where each $X_i$ has an arbitrary ID 
 distribution that need not be in the  same family as the other $X_{j\ne i}$,
 {\it i.e.}\ the individual distributions need not all be gamma or all fractional gamma with different parameter choices.
 They referred to the associated  normalised distribution generically  as the class of distributions on the 
  simplex. 
 
There is ample room for further exploration of multivariate alternatives to Dirichlet building upon the classes of densities 
explored in this paper. In addition to that, there is still much to explore in the world of univariate ID and GGC densities and the form of 
LKCD representation that they induce which, in turn, can lead to further insights.

\bibliography{Citations_2021}{}
\bibliographystyle{plain} 

\end{document}